\documentclass[12pt]{amsart}
\usepackage{amsmath,amssymb,latexsym}
\usepackage{comment,nicefrac}
\usepackage{tikz,pgfplots}
\usepackage{hyperref}
\usetikzlibrary{arrows,snakes,backgrounds}
\theoremstyle{plain}             
\newtheorem{theorem}{Theorem}[section]
\newtheorem{lemma}[theorem]{Lemma}

\newtheorem{example}[theorem]{Example}

\usepackage[normalem]{ulem}
\usepackage{dsfont}
\newcommand{\isdef}{\mathrel{\mathrel{\mathop:}=}}

\newcommand{\red}[1]{\textcolor{red}{#1}}

\newcommand{\spn}{\operatorname{span}}

\newcommand{\ALEoV}[2]{\mathbf{T}(#1,\mathbf{#2})}
\newcommand{\grey}[1]{\textcolor{black!70!white}{#1}}

\begin{document}
\title[An inverse transient Stefan problem] 
{An inverse problem for the one-phase Stefan problem 
with varying melting temperature}
\author{Marc Dambrine}
\address{Marc Dambrine,
Universit\'e de Pau et des Pays de l'Adour,
IPRA-LMA, UMR CNRS 5142, 
Avenue de l'universit\'e, 64000 Pau, France}
\email{marc.dambrine@univ-pau.fr}
\author{Helmut Harbrecht}
\address{Helmut Harbrecht,
Departement Mathematik und Informatik,
Universit\"at Basel,
Spiegelgasse 1, 4051 Basel, Schweiz}
\email{helmut.harbrecht@unibas.ch}

\maketitle
\begin{abstract}
The present article is dedicated to the forward and
backward solution of a transient one-phase Stefan problem.
In the forward problem, we compute the evolution of the initial 
domain for a Stefan problem where the melting temperature
varies over time. This occurs in practice, for example, when 
the pressure in the external space changes in time. In the 
corresponding backward problem, we then reconstruct the 
time-dependent melting temperature from the knowledge 
of the evolving geometry. We develop respective numerical 
algorithms using a moving mesh finite element method and 
provide numerical simulations.
\end{abstract}

\section{Introduction}
\subsection{Stefan problem}
The classic Stefan problem dates back to Josef Stefan 
in 1889, who studied ice formation in the polar seas, see 
\cite{Stefan1889}. He discovered that this growth problem 
is related to heat problems in which not only the temperature, 
which is the solution of a partial differential equation, is unknown
but also the position of the surface moves over time and is 
therefore part of the problem.

This particular problem belongs to the class of problems 
with moving boundaries. In general, such problems contain 
time-dependent boundaries that are unknown and depend 
on temporal and spatial variables. Problems with moving 
boundaries are also referred to as Stefan problems. This 
contrasts with problems with free boundary. The latter 
also contain boundaries that are unknown in advance, 
but these boundaries are in a steady state and are 
therefore not time-dependent, cf.~\cite{Crank1984}. 

There is a wealth of literature on the Stefan problem, see 
\cite{Friedman1975,Gupta2003,Meirmanov1992,Rubinstein1971, 
Tarzia1988} for example and the references contained therein. 
This literature deals mainly with the analysis of the Stefan 
problem. Stefan problems arise, for example, in the
modelling of phase transitions, chemical reactions, fluid 
flows in porous media or the melting of ice, compare 
\cite{Crank1984}.

When solving a Stefan problem numerically, one encounters the 
problem of handling the moving boundary. In \cite{Crank1984}, 
various numerical methods are explained, compare also
\cite{Furz,KBO} and the references therein: front tracking 
methods, front fixing methods, and fixed domain methods.
Front tracking methods calculate the position of the moving 
boundary at each time step. Front fixing methods attempt to 
fix the front by choosing a suitable spatial coordinate system. 
In fixed domain methods, the problem is reformulated, for 
example using the enthalpy method, so that the position 
of the boundary appears as a feature of the solution, see 
\cite{Crank1984}. Note that also shape optimization can 
be used to compute directly the evolving surface of the 
Stefan problem, compare \cite{Bruegger,Lujano}.

\subsection{Motivation and background}
As part of a recent collaboration with physicists, we looked in \cite{SDHDB}
at the growth and melting of methane hydrates. The experiment involved 
changing the pressure inside a capillary tube in which such a hydrate 
crystal grows in the presence of two components: water and methane gas 
dissolved in it. The growth or melting of the crystal was filmed and the rate 
of gas dissolved in the liquid part was measured by Raman spectrometry.  

The conclusion of the study \cite{SDHDB} is that the dynamics of the 
boundary and the gas concentrations measured are consistent with a 
two-phase Stefan model, in which the melting temperature at the interface 
between the hydrate and the aqueous phase is not the melting temperature 
of the hydrate single crystal. To reach this conclusion, we used, in the 
one-dimensional framework of the capillary experiments, analytical 
formulas giving the evolution of the front as a function of this temperature 
at the melting front. The experimental observation focused on the dynamics 
of the front, and an explicit formula was used to determine the temperature 
at the front.

In this article, we ask the following general question, which generalizes 
the previous situation in a context  where we do not have analytical formulas: 
is it possible to determine the melting temperature from observations of 
the evolution of the free surface? This melting temperature depends in 
particular on pressure. To limit the number of experiments required, the 
idea would be to vary this pressure slowly in the experimental setup. To 
do this, we must consider a transient version of Stefan's problem. 

The specific question we are considering in this article is therefore the 
following: if we observe the evolution of a free boundary evolving according 
to a Stefan model in which the temperature at the interface between phases 
varies with time, is it possible to reconstruct the melting temperature as a 
function of time?

\subsection{Our contribution}
Answering the above question in a general context would be ambitious 
and is not the objective of this article. We wish to conduct a preliminary 
study and respond from a numerical perspective to the feasibility of this 
problem in a proof-of-concept spirit. That is why we will place ourselves 
in the simplest non-trivial framework possible: the evolution of a star-shaped 
domain in two dimensions.  

For the direct problem, we consider a single-phase Stefan problem where 
the temperature imposed at the boundary is given and depends on time. 
The inverse problem is then as follows: suppose we know the evolution of 
a domain governed by a Stefan model with a temperature that varies over 
time. Can we reconstruct the melting temperature? 

\subsection{Content}
This article is organized as follows. First, in Section~\ref{sec:forward},
we introduce notation and specify the transient Stefan problem that 
we are considering. The numerical method that we propose for solving 
the problem under consideration is presented on Section~\ref{sct:forward}. 
Specifically, we employ a finite element method to solve the heat equation
in combination with a front tracking method. Section~\ref{sct:inverse} is 
then concerned with the inverse problem. We discuss necessary conditions 
for identifiability and propose some related positive and negative examples.
We then present numerical experiments that demonstrate that monitoring 
of the free boundary enables the determination of the varying temperature 
at the front. Finally, we state concluding remarks in Section~\ref{sct:conclusio}. 

\section{A nonstationary Stefan problem}\label{sec:forward}
\subsection{Classical one-phase Stefan problem}
\begin{figure}[htb]
\begin{center}
\begin{tikzpicture}[scale=0.6]
\filldraw[fill=black!20!white, densely dashed] (16,0) circle (2cm and 0.5cm);
\filldraw[fill=red!20!white, densely dashed] (15.85,3) circle (1.25cm and 0.37cm);
\draw[bend left=20] (14,6) to node [auto] {} (14,0);
\draw (14,0) arc (180:360:2cm and 0.5cm);
\draw[bend left=30] (18,0) to node [auto] {} (18,6);
\draw (18,6) ++ (-2,0) circle (2cm and 0.5cm);
\draw[densely dashed] (14,0) arc (180:0:2cm and 0.5cm);
\draw[densely dashed, purple] (15.85,3) circle (1.25cm and 0.37cm);

\node at (16,-2){$Q_T = \bigcup_{0 < t < T} \big( \{ t \} \times \red{\Omega_t} \big)$};
\node at (16,0){$\grey{\Omega_0}$};
\node at (16,3){$\textcolor{red}{\Omega_t}$};
\node at (14,3){$\textcolor{purple}{\Gamma_t}$};

\draw[black, thick,->] (10,0,0) -- (11.5,0,0) node[anchor=south]{$x_1$};
\draw[black, thick,->] (10,0,0) -- (10,1.5,0) node[anchor=north west]{$t$};
\draw[black, thick,->] (10,0,0) -- (10,0,1.5) node[anchor=north]{$x_2$}; 
\end{tikzpicture}
\end{center}
\caption{Geometric setup of the Stefan problem: The space-time
tube $Q_T$ is given by the time-evolution of the initial domain 
$\Omega_0$. The cut area at time $t$ is the domain $\Omega_t$
with boundary $\Gamma_t = \partial\Omega_t$.}
\label{fig.setup_stefan}
\end{figure}
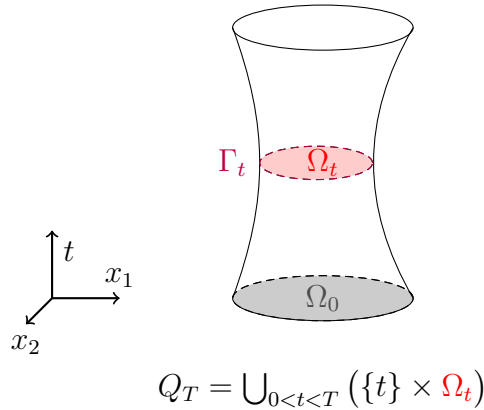

Let us consider the classical one-phase Stefan problem as described 
in \cite{Hadzic2014}. This specific Stefan problem models the evolution 
of the solid-liquid phase interface. Thus, for every point of time $t \in [0,T]$, 
we have a time-dependent spatial domain which we denote by $\Omega_t 
\subset \mathbb{R}^d$, $d \geq 2$. This spatial domain has a time-dependent 
spatial boundary $\Gamma_t \isdef \partial \Omega_t$. By setting 
\begin{equation}\label{eq:space-time tube}
Q_T = \bigcup_{0<t<T} \big( \lbrace t \rbrace \times \Omega_t \big), 
\end{equation}
we obtain the space-time non-cylindrical domain (also called 
\emph{space-time tube}) with lateral boundary
\[ 
\Sigma_T = \bigcup_{0<t<T} \big( \lbrace t \rbrace \times \Gamma_t \big).
\]
The setup is illustrated in Figure \ref{fig.setup_stefan} for two 
spatial dimensions plus the temporal dimension. 

For the formulation of the one-phase Stefan problem, we follow 
\cite{Hadzic2014, Hadzic2016}. The temperature $u(t,{\bf x})$ of 
the liquid in $\Omega_t$ is thus described by the partial differential 
equation
\begin{alignat}{3}
\label{eq.Stefan_PDE}
\partial_t u - \Delta u & = 0 \ \ && \text{ in } \Omega_t, \\
\label{eq.Stefan_cond}
\langle {\bf V}, {\bf n} \rangle & = - \frac{\partial u }{\partial {\bf n}} \ \ && \text{ on } \Gamma_t,\\
\label{eq.Stefan_boundary_cond}
u & = u_m \ \ && \text{ on } \Gamma_t, \\
\label{eq.Stefan_initial_cond}
u(0, \cdot ) & = u_0 \ \ && \text{ in } \Omega_0 = \Omega. 
\end{alignat}
We assume throughout this article, that the Dirichlet
data match the initial condition, meaning that $u(0,{\bf x})
= u_m$ for all ${\bf x}\in\Gamma_0$.

The domain $\Omega$ in \eqref{eq.Stefan_initial_cond} is the 
initial shape of the liquid phase while condition \eqref{eq.Stefan_cond} 
is called the Stefan condition \cite{Hadzic2016}. It comes from the 
movement of the phase interface, see \cite[pg.~387]{Visintin2008}. The 
Stefan condition expresses that the normal velocity $\langle {\bf V}, {\bf n} \rangle$ 
of the surface $\Gamma _t$ equals minus the normal derivative of $u$ at the boundary. 
We prescribe the initial position of the interface and the initial temperature distribution 
to make the problem meaningful. From this Stefan problem, we can see that the liquid 
freezes at zero temperature, cf.~\cite{Hadzic2014}. Notice that the one-phase Stefan problem 
is actually also a two-phase Stefan problem, but the temperature is only unknown in 
one region, while it is equal to $u_m$ in the other region, compare \cite{Visintin2008}.

The domain $\Omega_t$, thus the region which contains the liquid phase, 
is characterized by $\lbrace {\bf x} \in \mathbb{R}^d: \; u(t, {\bf x}) > u_m \rbrace$ 
if we have $u_0 > u_m$ in $\Omega$. Therefore, $u$ can be interpreted as 
a level set function. Due to \eqref{eq.Stefan_PDE}, the parabolic Hopf 
lemma (see e.g.~\cite{Friedman1958} for some remarks) implies 
$\partial u/\partial {\bf n} < 0$ on $\Gamma_t$ for $t>0$. Therefore, 
we obtain the so-called \emph{Rayleigh-Taylor sign condition}
\[ 
- \frac{\partial u_0}{\partial {\bf n}} \geq \lambda > 0 \text{ on } \Gamma_0, 
\]
which ensures the nondegeneracy in accordance with \cite{Hadzic2014}.
Vice versa, when $u_0 < u_m$ in $\Omega$, we find $\partial u/\partial {\bf n} < 0$ 
on $\Gamma_t$ for $t>0$ and hence
\[ 
- \frac{\partial u_0}{\partial {\bf n}} \leq \lambda < 0 \text{ on } \Gamma_0.
\]

\subsection{Varying the melting temperature}
The melting temperature in the classical Stefan problem 
\eqref{eq.Stefan_PDE}--\eqref{eq.Stefan_initial_cond} is fixed.
In contrast, we consider in the present article the situation that 
the melting temperature $u_m = u_m(t)$ varies over time and 
is thus time-dependent. In practice, this can be achieved by 
changing the pressure in the external space during time. 

In order to solve the Stefan problem 
\eqref{eq.Stefan_PDE}--\eqref{eq.Stefan_initial_cond} with 
time-dependent melting temperature, we make the ansatz 
$v=u-u_m$. In view of $\Delta u_m = 0$, this amounts to
\begin{alignat}{3}
\label{eq.Stefan_PDE'}
\partial_t v - \Delta v & = -\dot{u}_m \ \ && \text{ in } \Omega_t, \\
\label{eq.Stefan_cond'}
\langle {\bf V}, {\bf n} \rangle & = - \frac{\partial v}{\partial {\bf n}} \ \ && \text{ on } \Gamma_t,\\
\label{eq.Stefan_boundary_cond'}
v & = 0 \ \ && \text{ on } \Gamma_t, \\
\label{eq.Stefan_initial_cond'}
v(0, \cdot ) & = v_0\ \ && \text{ in } \Omega_0 = \Omega,
\end{alignat}
where we set $v_0\isdef u_0-u_m(0)$. 

In the following, we will solve the 
transformed problem \eqref{eq.Stefan_PDE'}--\eqref{eq.Stefan_initial_cond'} 
rather than the original problem \eqref{eq.Stefan_PDE}--\eqref{eq.Stefan_initial_cond}.
The solution $u$ of the latter is then obtained in accordance with $u = v+u_m$.

\subsection{Generation of the space-time tube}
In order to generate a space-time tube, we can adopt two different 
points of view. For both of them, let us assume that we have a spatial 
domain $\Omega_{\text{ref}}$, which serves as reference domain. One 
can generate a space-time tube $Q_T$ by mapping this domain 
$\Omega_{\text{ref}}$ to a spatial domain $\Omega_t$ for every 
point of time $t$, see also Figure~\ref{fig.setup_stefan}. 

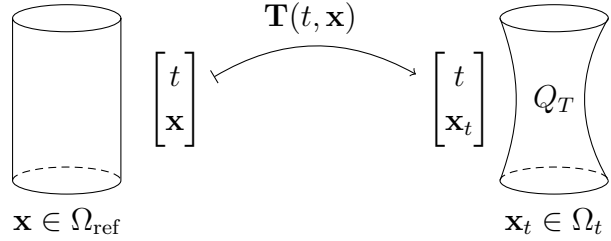
\begin{figure}[htb]
	\centering
	\begin{tikzpicture}[scale=0.36]
		\draw (-2,6) -- (-2,0) arc (180:360:2cm and 0.5cm) -- (2,6) ++ (-2,0) circle (2cm and 0.5cm);
		\draw[densely dashed] (-2,0) arc (180:0:2cm and 0.5cm);
		
		\draw[bend left=20] (16,6) to node [auto] {} (16,0);
		\draw (16,0) arc (180:360:2cm and 0.5cm);
		\draw[bend left=30] (20,0) to node [auto] {} (20,6);
		\draw (20,6) ++ (-2,0) circle (2cm and 0.5cm);
		\draw[densely dashed] (16,0) arc (180:0:2cm and 0.5cm);
		
		\node (A) at (4,3) {$\begin{bmatrix} t \\ {\bf x} \end{bmatrix}$};
		\node (G) at (14.5,3) {$\begin{bmatrix} t \\ {\bf x}_t \end{bmatrix}$};
		\node (B) at (18,3) {$Q_T$};
		\node (B) at (18,-1.5) {${\bf x}_t \in \Omega_t$};
		\node (C) at (0,-1.5) {${\bf x} \in \Omega_{\text{ref}}$};
		
		\draw[|->,  bend left=30] (A) to (G);
		\node (E) at (9,6) {${\bf T} (t, {\bf x}) $};
	\end{tikzpicture}
	\caption{Generation of the space-time tube by the mapping 
	${\bf T}(t, \cdot)$ induced by the velocity field ${\bf V}$. In 
	particular, the relation $\Omega_t = {\bf T}(t,\Omega_{\text{ref}})$ 
	always applies.}
	\label{fig.generate_tube_velocity}
\end{figure}

On the one hand, this can be done by considering a velocity field 
${\bf V}$ and associate to it the solution $\ALEoV{t}{\cdot} : {\bf x}
\mapsto {\bf x}_t = \ALEoV{t}{x}$ of the differential equation 
\begin{equation} \label{eq.speed_method}
	\begin{aligned}
		\frac{\partial }{\partial t} \ALEoV{t}{x}  
		&= {\bf V}\big(t, \ALEoV{t}{x}\big) && \ \ \text{in }  (0,T) \times \Omega_{\text{ref}},  \\
		\ALEoV{0}{x} & = {\bf T}_0({\bf x}) && \ \ \text{in } \Omega_{\text{ref}},
	\end{aligned}
\end{equation}
see \cite[pg.~6]{Zolesio1979} or \cite{Moubachir2006} for the 
details. The map ${\bf T}(t,\cdot)$ thus describes the pathline of 
an individual particle being exposed to the velocity field ${\bf V}$. 
By setting $\Omega_t = {\bf T} (t, \Omega_{\text{ref}})$, we 
finally generate the space-time tube $Q_T$, see Figure 
\ref{fig.generate_tube_velocity} for an illustration.

\subsection{Solution algorithm for the transient Stefan problem}
In order to derive an algorithm for the transient Stefan problem 
\eqref{eq.Stefan_PDE'}--\eqref{eq.Stefan_initial_cond'}, we choose 
a small time step $\Delta t$ and set $t_k\isdef k\Delta t$ for 
$k\in\mathbb{N}_0$. We then aim at approximating the domains
$\Omega_k\isdef\Omega_{t_k}$ and hence the mappings
${\bf T}_k\isdef {\bf T}(t_k,\cdot)$ for all $k\in\mathbb{N}$. 

To this end, we solve successively the heat equation for each 
time interval $t\in [t_k,t_{k+1}]$ on the fixed domain $\Omega_k$.
This means, we consider the time slice $[t_k,t_{k+1}]\times\Omega_k$
and solve the initial boundary value problem
\begin{alignat}{3}
\label{eq.Stefan_PDE"}
\partial_t \widetilde{v} - \Delta \widetilde{v} & = -\dot{u}_m \ \ && \text{ in } \Omega_k, \\
\label{eq.Stefan_boundary_cond"}
\widetilde{v} & = 0 \ \ && \text{ on } \Gamma_k, \\
\label{eq.Stefan_initial_cond"}
\widetilde{v}(t, \cdot ) &= v_k\ \ && \text{ in } \Omega_k,
\end{alignat}
where $v_k$ is the solution from the previous time step
at time $t_k$. After having the solution $\widetilde{v}$ at
hand, we can use the Stefan condition \eqref{eq.Stefan_cond'} 
to update the boundary of the (approximate) space-time tube 
in accordance with 
\begin{equation}
\label{eq.Stefan_cond"}
  \Gamma_{k+1}\isdef\bigg\{{\bf x}-\Delta t
  \frac{\partial \widetilde{v}}{\partial {\bf n}}(t_{k+1},{\bf x}){\bf n}({\bf x}):
  {\bf x}\in\Gamma_k\bigg\}.
\end{equation}
This defines the new domain $\Omega_{k+1}$ and thus the
map ${\bf T}_{k+1}$. The new initial data $v_{k+1}$, entering
the computation of $\widetilde{v}$ in the next time slice $[t_{k+1},
t_{k+2}]\times \Omega_{k+1}$, are obtained by interpolating the solution 
$\widetilde{v}$ of \eqref{eq.Stefan_PDE"}--\eqref{eq.Stefan_initial_cond"} 
at time $t = t_{k+1}$ onto the updated domain $\Omega_{k+1}$. Hence, 
we especially have the identity
\[
  v_{k+1}({\bf x})\isdef \widetilde{v}(t_{k+1},{\bf x})
  \quad\text{for all ${\bf x}\in\Omega_k\cap\Omega_{k+1}$}.
\]

Having the new initial data $v_{k+1}$ and the updated 
domain $\Omega_{k+1}$ at hand, we can repeat the above 
procedure to compute $\widetilde{v}$ on the next time slice 
$[t_{k+1},t_{k+2}]\times\Omega_{k+1}$. Since we track 
the boundary of the sought space-time cylinder $Q_T$, 
associated with the transient Stefan problem 
\eqref{eq.Stefan_PDE}--\eqref{eq.Stefan_initial_cond},
the proposed method belongs to the class of front 
tracking methods.

\section{Numerical simulation of the forward solution}
\label{sct:forward}
\subsection{Discretization of the space-time tube}
The key ingredient of any numerical 
method to solve the transient Stefan problem 
\eqref{eq.Stefan_PDE}--\eqref{eq.Stefan_initial_cond}
is the discretization of the sought space-time tube $Q_T$.
As we approximate it by a piecewise constant function in
time through the approximations of the time slices $[t_k,
t_{k+1}]\times\Omega_k$, $k\in\mathbb{N}_0$, we only 
need to discretize the domain $\Omega_k$, which is
considered in the following.

Under the assumption that the domain $\Omega_k$ 
is star-shaped, we can represent its boundary by means 
of a Fourier series. We hence can make the ansatz
\begin{equation}\label{eq:boundary}
  \Gamma_k\isdef\Bigg\{\bigg(a_{k,0} 
  	+ \sum_{\ell=1}^M \big[a_{k,\ell}\cos(\ell\varphi) 
  	+ a_{k,-\ell}\sin(\ell\varphi)\big]\bigg)
	\begin{bmatrix}\cos(\varphi)\\\sin(\varphi)\end{bmatrix}: 
	\varphi\in [0,2\pi]\Bigg\},
\end{equation}
which identifies the domain $\Omega_k$ (respectively its boundary 
$\Gamma_k$) with the $2M+1$ Fourier coefficients
\[
{\bf a}_k = [a_{k,-M},\ldots,a_{k,-1},a_{k,0},a_{k,1},\ldots,a_{k,M}]
\in\mathbb{R}^{2M+1}
\]
of the representation of its boundary $\Gamma_k$. 

If we define the reference domain $\Omega_{\text{ref}}$ as
the unit disc, that is
\[
\Omega_{\text{ref}} := \{{\bf x}\in\mathbb{R}^2:\|{\bf x}\|_2 < 1\}.
\]
we obtain for the associated domain $\Omega_k$ the 
mapping ${\bf T}_k$ by
\begin{equation}\label{eq:T_k}
  {\bf T}_k:\Omega_{\text{ref}}\to\Omega_k,\quad
  {\bf T}_k({\bf x}) = \bigg(a_{k,0} + \sum_{\ell=1}^M \big[a_{k,\ell}\cos\big(\ell\widehat{\bf x}\big) 
  	+ a_{k,-\ell}\sin\big(\ell\widehat{\bf x}\big)\big]\bigg){\bf x}.
\end{equation}
Herein, the expression
\[
 \widehat{\bf x}\isdef\frac{\bf x}{\|{\bf x}\|}
\]
denotes the radial direction, which is interpreted as the polar angle in 
the series expansion. 

\begin{figure}[hbt]
\begin{center}
\includegraphics[width=0.45\textwidth,trim={170 40 160 30},clip]{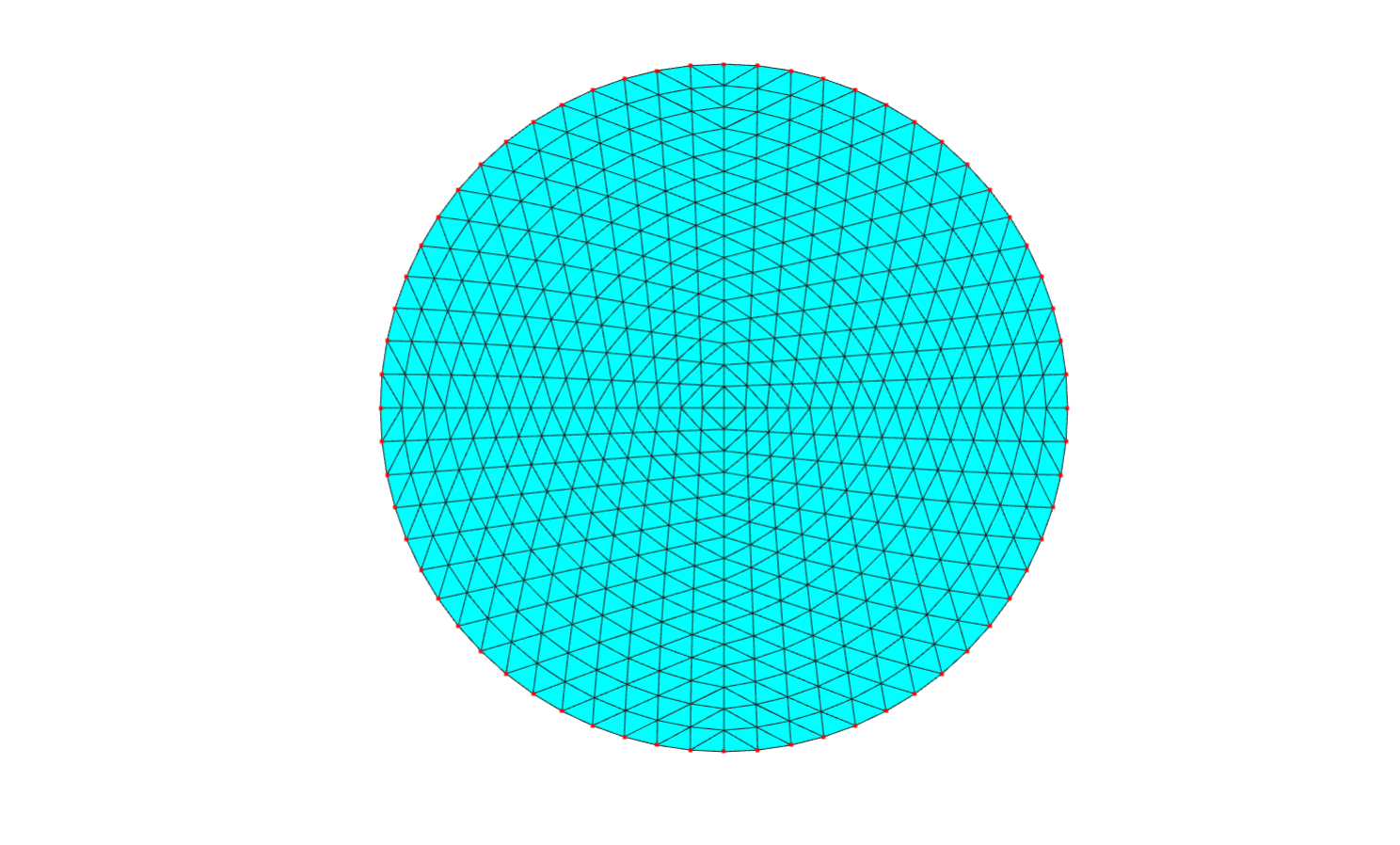}
\includegraphics[width=0.45\textwidth,trim={170 40 160 30},clip]{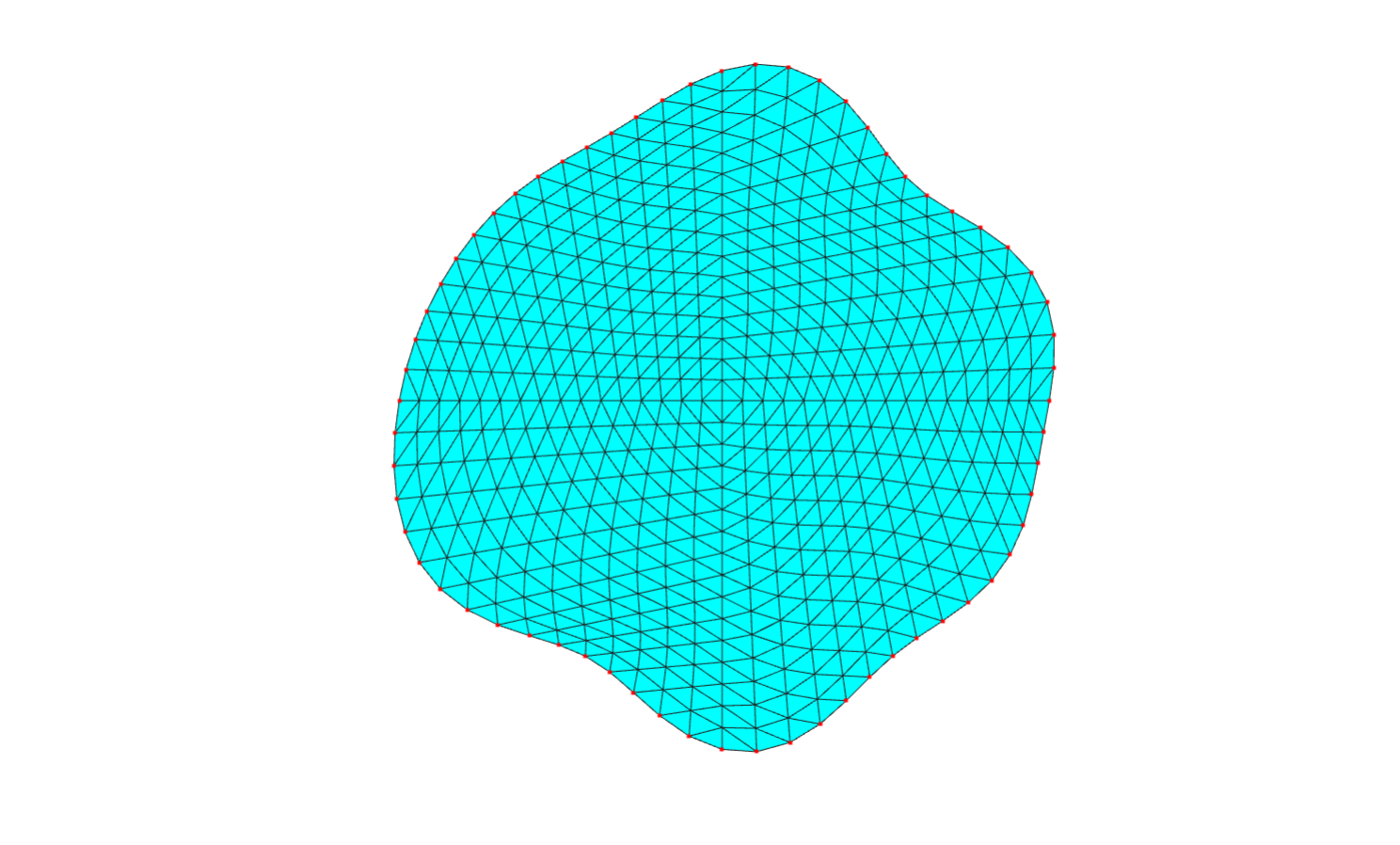}
\caption{\label{fig:mesh}
The finite element mesh on the circle and its mapped 
counterpart on the deformed domain $\Omega_k$.}
\end{center}
\end{figure}

We refer to Figure~\ref{fig:mesh} for an illustration of our construction 
of the mapping ${\bf T}_k$ including finite element meshes on 
$\Omega_{\text{ref}}$ and on the given domain $\Omega_k$. Note
that the finite element mesh on $\Omega_k$ is obtained by mapping 
the finite element mesh on $\Omega_{\text{ref}}$ by using ${\bf T}_k$.
Note that the boundary vertices of the mesh on the unit disc 
$\Omega_{\text{ref}}$ are equidistant which ensures that the boundary
vertices of the mesh on the mapped domain $\Omega_k$ are always 
equidistant with respect to the polar angle. This enables the application
of the fast Fourier transform \cite{FFT} to speed up the computations
in the following.

\subsection{Discretization of the heat equation}
In order to solve the heat equation 
\eqref{eq.Stefan_PDE"}--\eqref{eq.Stefan_initial_cond"}
in a given time slice $[t_k,t_{k+1}]\times\Omega_k$, we
apply the finite element method. To this end, we compute a 
mesh by mapping a given quasi-uniform mesh of the unit circle 
$\Omega_{\text{ref}}$ to the actual domain $\Omega_k$ by 
application of the mapping ${\bf T}_k$ from \eqref{eq:T_k}. 
Note that this mesh mapping method ensures that we always 
have the same mesh with vertices that are transported by ${\bf T}_k$. 

For sake of simplicity, we assume that the melting temperature
$u_m(t)$ is piecewise linear in each time interval $[t_k,t_{k+1}]$.
Then, given the domain $\Omega_k$ in the $(k+1)$-st time step, 
we consider a quasi-uniform mesh (computed as described above)
and the associated finite element space $V_k = \spn\{\varphi_i: 
i=1,\ldots,N\}$, consisting of $N$ piecewise linear Lagrangian 
finite element basis functions $\{\varphi_i\}_i$. We define the 
finite element stiffness and mass matrices
\[
  {\bf A}_k = \big[(\nabla\varphi_i,\nabla\varphi_j)_{L^2(\Omega_k)}\big]_{i,j=1}^N,\quad
  {\bf M}_k = \big[(\varphi_i,\varphi_j)_{L^2(\Omega_k)}\big]_{i,j=1}^N,
\]
as well as the right-hand side
\[
  {\bf f}_k = \big[(1,\varphi_j)_{L^2(\Omega_k)}\big]_{i=1}^N,
\]
compare \cite{Braess,Brenner/Scott} for example. Then, 
the finite element function 
\[
v_{k+1} = \sum_{i=1}^N v_{k+1,i}\varphi_i\in V_k
\] 
is given by using the trapezoidal rule in time as the 
solution of the following linear system of equations
\begin{equation}\label{eq:LSE}
  \bigg({\bf M}_k+\frac{\Delta t}{2}{\bf A}_k\bigg){\bf v}_{k+1}
  = \bigg({\bf M}_k-\frac{\Delta t}{2}{\bf A}_k\bigg){\bf v}_k-\Delta t \alpha_k {\bf f}_k.
\end{equation}
Herein, $\alpha_k$ denotes the derivative of the melting 
temperature $u_m(t)$ in the interval $[t_k,t_{k+1}]$. We
emphasize $\alpha_k$ can be replaced by an appropriate 
constant approximation of $\dot{u}_m$ if the melting temperature 
$u_m(t)$ is not linear in the time interval $[t_k,t_{k+1}]$.
In the context of finite element methods, the discretization
underlying the system \eqref{eq:LSE} is also known as 
Crank-Nicolson method, introduced in \cite{CK} and 
being second order accurate in space and time.

\subsection{Updating the domain}
Having the finite element solution $v_{k+1}$ for the time
step $t_{k+1}$ at hand, we can compute the update of the
domain $\Omega_k$. To this end, we note that the gradient
of $v_{k+1}$ at the boundary $\Gamma_k$ is given by
\[
\nabla v_{k+1}({\bf x}) = \frac{\partial v_{k+1}}{\partial{\bf n}}({\bf x}){\bf n}({\bf x}),
\]
since the tangential derivative $(\partial v_{k+1}/\partial{\bf t})({\bf x}) 
= 0$ vanishes for all ${\bf x}\in\Gamma_k$ due to the homogenous 
Dirichlet data. Therefore, we can update the boundary in radial direction 
by using the modified Stefan condition
\[
  -\langle {\bf V},\widehat{\bf x}\rangle 
  = \langle\nabla v_{k+1},\widehat{\bf x}\rangle 
  = \frac{\partial v_{k+1}}{\partial{\bf n}}\langle{\bf n},\widehat{\bf x}\rangle
  \quad\text{on $\Gamma_k$}.
\]
We hence have to update all points on the boundary in the 
radial direction by means of the formula
\begin{equation}\label{eq:update}
  {\bf x}_{k+1}\isdef {\bf x}-\Delta t \frac{\partial v_{k+1}}{\partial{\bf n}}
  	\langle{\bf n},\widehat{\bf x}\rangle\widehat{\bf x}\quad\text{on $\Gamma_k$}.
\end{equation}

In our numerical method, we realize the update \eqref{eq:update}
as follows. Under the assumption that the vertices of the mesh on 
the boundary correspond to an equidistant subdivision of the radial 
function with respect to the polar angle, we can apply the fast Fourier 
transform to compute the update of the coefficients in the Fourier
series. However, there are typically more vertices on the boundary 
than Fourier coefficients in the representation \eqref{eq:boundary}.
Therefore, we always truncate the update to its $2M+1$ leading 
coefficients which corresponds to a least squares fit of the new 
parametrization using the $2M+1$ Fourier coefficients.

In a final step, the finite element function $v_{k+1}$, given on the 
old domain $\Omega_k$, is interpolated onto the updated domain 
$\Omega_{k+1}$ by means of kernel interpolation. Kernel interpolation 
is a very efficient approach for scattered data interpolation, see 
\cite{Fasshauer,Wendland} for example. In our implementation, 
we use the exponential kernel $\kappa({\bf x},{\bf y}) = 
\exp(-\|{\bf x}-{\bf y}\|_2)$ for the interpolation, which results 
in an interpolation error that is known to converge in 
the $L^2$-norm with order $h^{-2}$.

In order to explain the kernel interpolation method, let $\{{\bf y}_{k,i}\}_i$ 
denote the vertices of the finite element mesh on $\Omega_k$ and
let $v_{k+1}^{\text{old}}\in V_k$ be the finite element function 
$v_{k+1}$ which comes from the solution of \eqref{eq:LSE}.
Assuming that there are $L$ vertices, we define the $L\times L$ 
kernel matrices
\[
  {\bf K} = \big[\kappa({\bf y}_{k,i},{\bf y}_{k,j})\big]_{i,j=1}^L,
  \quad {\bf L}:=\big[\kappa({\bf y}_{k,i},{\bf y}_{k+1,j})\big]_{i,j=1}^L.
\]
Then, the finite element function $v_{k+1}^{\text{new}}\in V_{k+1}$ 
with respect to the new vertices $\{{\bf y}_{k+1,i}\}_i$ is given 
through the nodal values obtained from
\[
  {\bf v}^{\text{new}} = {\bf L}  {\bf K}^{-1} {\bf v}^{\text{old}},
  \quad \text{where}\ {\bf v}^{\text{old}} = \big[v_{k+1}^{\text{old}}({\bf y}_{k,i})\big]_{i=1}^L.
\]
This means that there holds $v_{k+1}^{\text{new}}({\bf y}_{k+1,i}) 
= [{\bf v}^{\text{new}}]_i$ for all new vertices ${\bf y}_{k+1,i}$, 
$i=1,\ldots,L$. We refer to \cite{Fasshauer,Wendland} for 
further details.

\subsection{Numerical results of the forward simulation}
\begin{figure}[hbt]
\begin{center}
\includegraphics[width=0.6\textwidth,trim={60 30 60 30},clip]{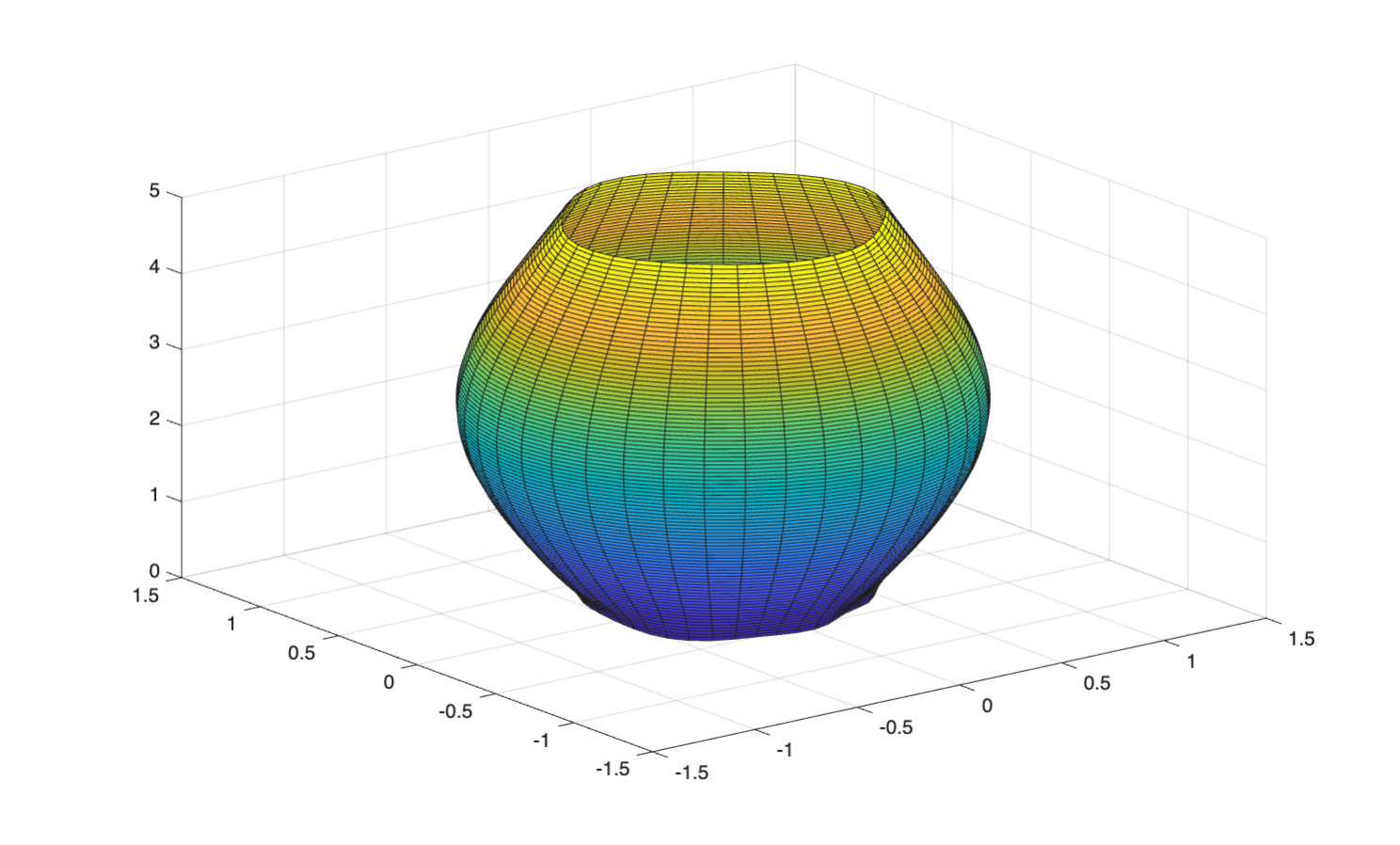}
\caption{\label{fig:example1}
The space-time tube in case of the melting temperature $u_m(t) 
= \nicefrac{1}{20}(t-\nicefrac{5}{2})^2$.}
\end{center}
\end{figure}

\begin{figure}[hbt]
\begin{center}
\includegraphics[width=0.6\textwidth,trim={60 30 60 30},clip]{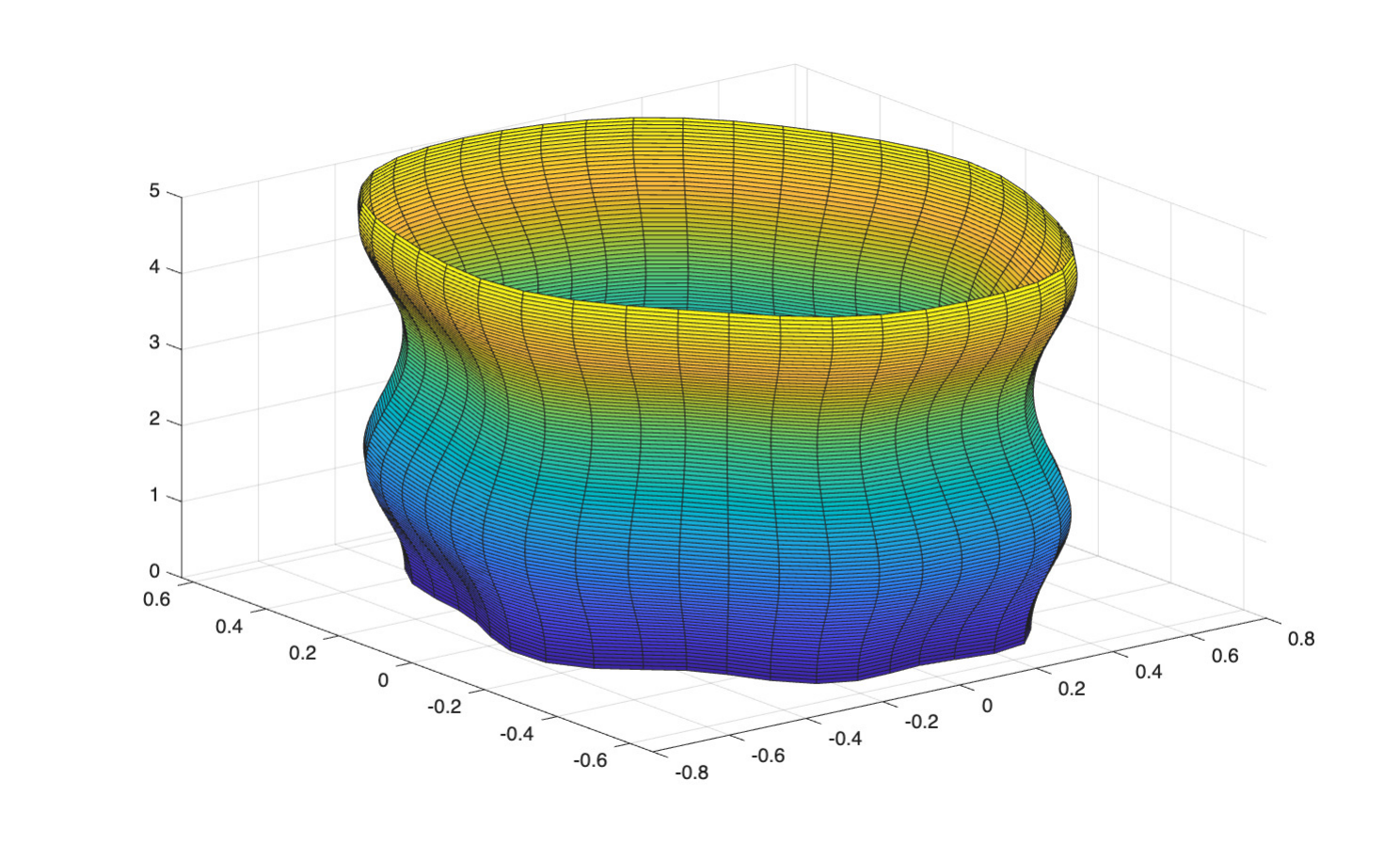}
\caption{\label{fig:example2}
The space-time tube in case of the melting temperature $u_m(t) 
= \nicefrac{1}{20}\big(\cos(2t)-1\big)$ on the right.}
\end{center}
\end{figure}

We shall present some numerical tests. We consider the time interval
$[0,5]$, that is $T=5$, and set $\Delta t = \nicefrac{5}{100}$, which 
results in 100 time steps. We apply $2M+1 = 29$ Fourier coefficients in 
the Fourier series of the radial function in the boundary representation 
\eqref{eq:boundary}. The finite element mesh we use consists of about
16000 finite elements, leading to about $N = 8000$ piecewise linear 
Lagrangian finite element basis functions. We prescribe the melting 
temperature $u_m(t) = \nicefrac{1}{20}(t-\nicefrac{5}{2})^2$ in the first 
example while it is set as $u_m(t) = \nicefrac{1}{20}\big(\cos(2t)-1\big)$ 
in the second example. The initial domain is chosen randomly in both 
simulations and the initial temperature has been set to be constant, 
which means $u_0\isdef u_m(0)$. The simulations run quite fast and 
the computed space-time tubes are found in Figure~\ref{fig:example1} 
for the first example and in Figure~\ref{fig:example2} for the second 
example.

In the first example, the domain grows for $t\in[0,\nicefrac{5}{2}]$,
becoming more and more a circle. Then, it shrinks again for $t\in
[\nicefrac{5}{2},5]$, compare Figure~\ref{fig:example1}. This is due 
to the increase of the melting temperature in the time interval 
$[0,\nicefrac{5}{2}]$ and the decrease of the melting temperature 
in the time interval $[\nicefrac{5}{2},5]$. We observe especially that 
the final boundary is close to the boundary of the initial domain. 

In the second example, the shape oscillates around the 
initial shape, while the shape of the surface becomes smoother
and smoother over time, compare Figure~\ref{fig:example2}. 
Indeed, the evolving boundary tends towards a disc. The reason for this 
is that the initial domain does not grow so much in comparison 
to the first example, so that it changes its shape only slightly.

\section{Inverse reconstruction}
\label{sct:inverse}
\subsection{On necessary conditions for the identifiability}
We shall next consider the \emph{inverse problem}, 
namely the reconstruction of the melting temperature 
from measurements of the evolving surface. 

The melting temperature appears as the Dirichlet trace 
of the temperature on the surface of the space-time tube. 
If the space-time tube is known, the deformation velocity 
of its surface is also known. Thus, if the space-time tube in 
question corresponds to Stefan's dynamics, the inverse 
problem of reconstructing the melting temperature is easy 
to solve. We should, however, take a closer look on the
identifiability.

Consider a given domain tube $Q_T$, given as in 
\eqref{eq:space-time tube}, which has been generated by 
Stefan's dynamics. In a first step, it suffices to solve the 
Neumann problem for the heat equation on $Q_T$, with 
Neumann data corresponding to the (negative) normal 
speed of deformation of the domain tube:
\begin{alignat*}{3}
\partial_t u - \Delta u & = 0 \ \ && \text{ in } \Omega_t, \\
 \frac{\partial u }{\partial {\bf n}} &=-\langle {\bf V}, {\bf n} \rangle\ \ && \text{ on } \Gamma_t,\\
u(0, \cdot ) & = u_0 \ \ && \text{ in } \Omega_0 = \Omega. 
\end{alignat*}
This problem has already been studied, notably in the article
of Hofmann and Lewis \cite{HofmannLewis} who established 
the existence and uniqueness of a solution to this problem 
under certain assumptions of regularity of the tube (see 
also the article of Br\"ugger et al.~in \cite{Bruegger2022} for 
an integral equation resolution). Then, in a second step, we 
must take the Dirichlet trace, which is necessarily the 
function to be reconstructed.

This procedure can, of course, be applied to any sufficiently 
regular space-time tube $Q_T$, even if it does not correspond 
to a Stefan-type evolution. The criterion that distinguishes the 
two cases is as follows: the trace of the solution to the heat 
equation on the space-time tube with the Neumann condition 
must be a constant on the boundary of each section $\Omega_t$ 
of the tube. 

\begin{example}
Let us consider the specific case of an disk 
\[
\Omega_0 = \{{\bf x}:\|{\bf x}\|_2<R_0\}
\]
as initial domain. The invariance of the problem with respect 
to rotations around the center of the disk implies that Stefan's
dynamics generates a space-time tube $Q_T$ whose sections
$\Omega_t$ are also disks.
As a consequence, there are space-time tubes that do not 
correspond to Stefan's dynamics: Let $\Omega_t$ be the 
ellipse with center at the origin that satisfies the equation
\[
\frac{x^2}{(1+t)^2}+\frac{y^2}{(1+2t)^2}=1.
\]
Then, the initial section $\Omega_0$ of the associated 
space-time tube $Q_T$ is a disk, but the domains 
$\Omega_t$ are not for $t>0$. 
\end{example}

We next like to discuss whether the knowledge of the 
final domain $\Omega_T$ allows us to reconstruct the 
temperature trajectory. To this end, let us consider a regular 
function $t\mapsto R(t)$ with non-negative values and the tube 
\[
  Q_T = \{(t,{\bf x}): \|{\bf x}\|_2<R(t),\ t\in[0,T]\}. 
\]  
The solution $u$ to the heat equation
\begin{alignat*}{3}
\partial_t u - \Delta u & = 0 \ \ && \text{ in } \Omega_t = \{{\bf x}:\|{\bf x}\|_2 < R(t)\}, \\
 \frac{\partial u }{\partial {\bf n}} &=-R'(t)\ \ && \text{ on } \Gamma_t,\\
u(0, \cdot ) & = u_0 \ \ && \text{ in } \Omega_0 = \Omega. 
\end{alignat*}
is radial, and therefore the transient melting temperature is 
simply $t\mapsto u\big(t,{\bf x}(t)\big)$ with $\|{\bf x}(t)\|=R(t)$.
As a consequence, we get

\begin{lemma}
For any $R_T>0$ and for any final time $T$, there 
exists a temperature trajectory $u_m(t)$ that transforms, 
through the Stefan dynamics, the initial disk into a disk with 
radius $R_T$ at time $T$.
\end{lemma}

The answer to the question of unique identifiability from
knowledge of the final domain $\Omega_T$ is thus clearly 
no. Consider the disk of radius $R_0$ as initial domain and 
the disk of radius $R_T$ as the final one at time $T$. For
any $R>0$, one can construct a Stefan tube such that 
$\Omega_0=\{{\bf x}:\|{\bf x}\|_2<R_0\}$, $\Omega_{T/2}
=\{{\bf x}:\|{\bf x}\|_2<R\}$ and $\Omega_T = \{{\bf x}:
\|{\bf x}\|_2<T\}$ by the previous lemma.

Since it generally appears difficult to characterize the necessary 
condition theoretically in geometric terms on a space-time tube, 
we will instead conduct numerical simulations.

\subsection{Reconstruction of the melting temperature}
\label{scr:inverse}
To reconstruct of the melting temperature from measurements of the 
evolving surface, we proceed as follows. For each time step $k=0,1,2,\ldots$, 
we solve \eqref{eq.Stefan_PDE"}--\eqref{eq.Stefan_initial_cond"} by setting
$u_m'\isdef 1$, meaning that we assume that the change of the melting
temperature in the time-interval $[t_k,t_{k+1}]$ is piecewise linear.
Since the solution $\widetilde{v}$ of 
\eqref{eq.Stefan_PDE"}--\eqref{eq.Stefan_initial_cond"} depends 
linearly on the inhomogeneity $-\dot{u}_m = -1$ and so does the 
flux, we can determine the scaling factor $\alpha_k$ such that 
we match the measured boundary $\Gamma_{k+1}$ by the 
update rule \eqref{eq.Stefan_cond"} best possible, meaning 
that we search for $\alpha_k\in\mathbb{R}$ such that
\begin{equation}\label{eq:identification}
  \Gamma_{k+1}\overset{!}{=}\bigg\{{\bf x}-\alpha_k\Delta t
  \frac{\partial \widetilde{v}}{\partial {\bf n}}(t_{k+1},{\bf x}){\bf n}({\bf x}):
  {\bf x}\in\Gamma_k\bigg\}.
\end{equation}
In that way, we can determine successively the time derivative
$\dot{u}_m$ of the melting temperature, leading to
\[
  \dot{u}_m(t) = \alpha_k, \quad \text{for all $t\in [t_k,t_{k+1}]$
  and $k\in\mathbb{N}_0$}.
\]
Therefore, by integrating this piecewise approximation 
with respect to the time, we obtain the sought melting 
temperature $u_m$. 

The method to solve equation \eqref{eq:identification},
which corresponds to a shape identification problem, 
is discussed in the next subsection. Especially, we 
investigate there the sensitivity of this reconstruction
with respect to measurement errors by numerical tests.

\subsection{Setup and algorithm}
We now test our reconstruction algorithm proposed in
the previous subsection. To that end, we consider the same
two examples as in Section~\ref{sct:forward}, but now try to 
reconstruct the melting temperature from the observed 
space-time tubes. To this end, we add (uncorrelated) random 
noise to the Fourier coefficients of the desired space-time 
tubes. The noise levels $\delta$ under consideration are 
$\delta = 0.25\%$, $\delta = 0.5\%$, $\delta = 1.0\%$, 
and $\delta = 2.0\%$. Recall that the prescribed space-time 
tubes have been visualized in Figure~\ref{fig:example1}
and Figure~\ref{fig:example2}, respectively. Moreover,
to avoid an inverse crime, we use only $2M+1 = 15$ 
Fourier coefficients for the reconstruction of the evolving
boundary of the space-time tube and a coarser mesh for
the finite element solver. It consists of about 2000 finite
elements, which results in about $N = 1000$ piecewise 
linear Lagrangian finite element basis functions. 

\begin{figure}[hbt]
\begin{center}
\includegraphics[width=0.48\textwidth,trim={60 30 60 30},clip]{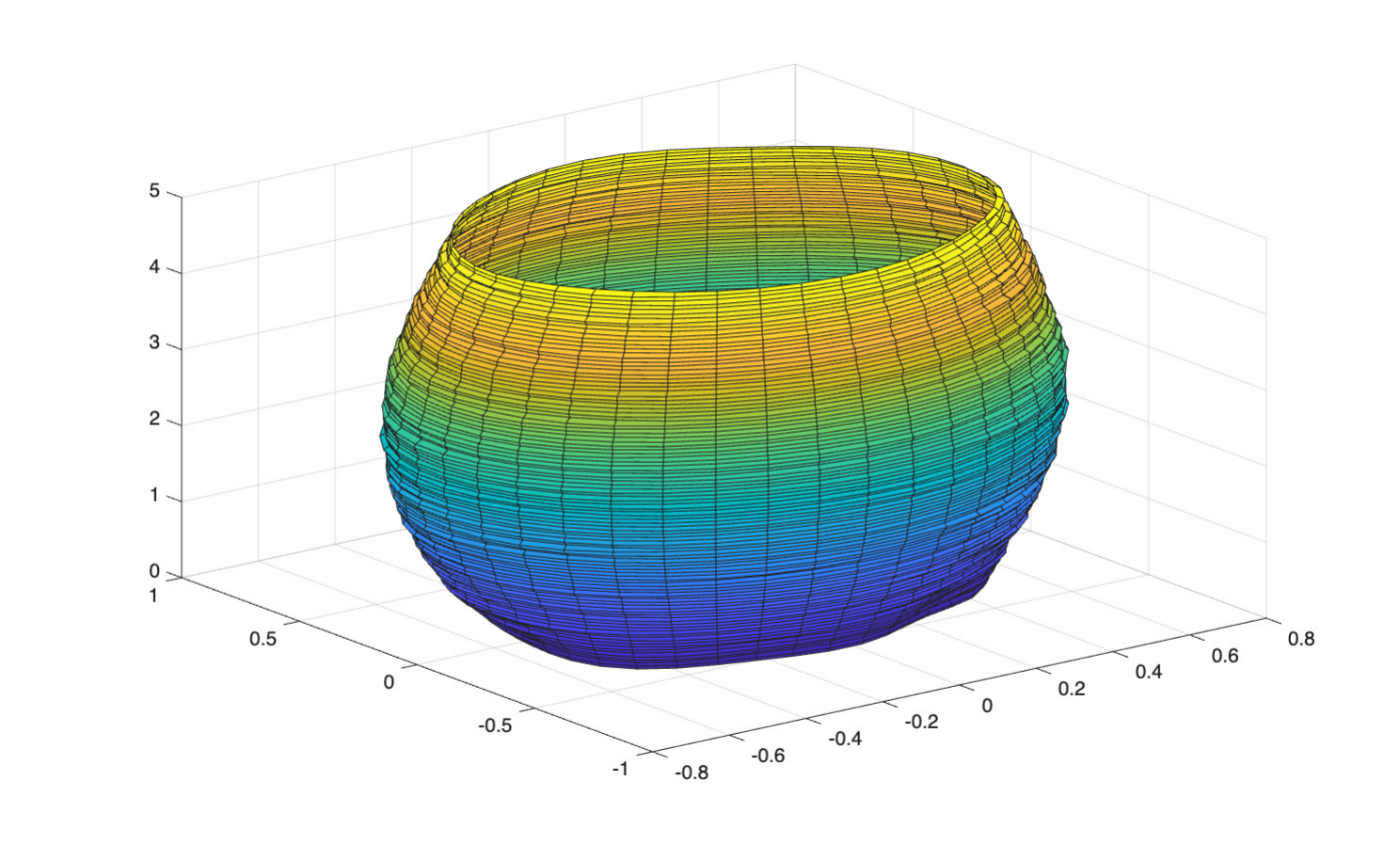}
\includegraphics[width=0.48\textwidth,trim={60 30 60 30},clip]{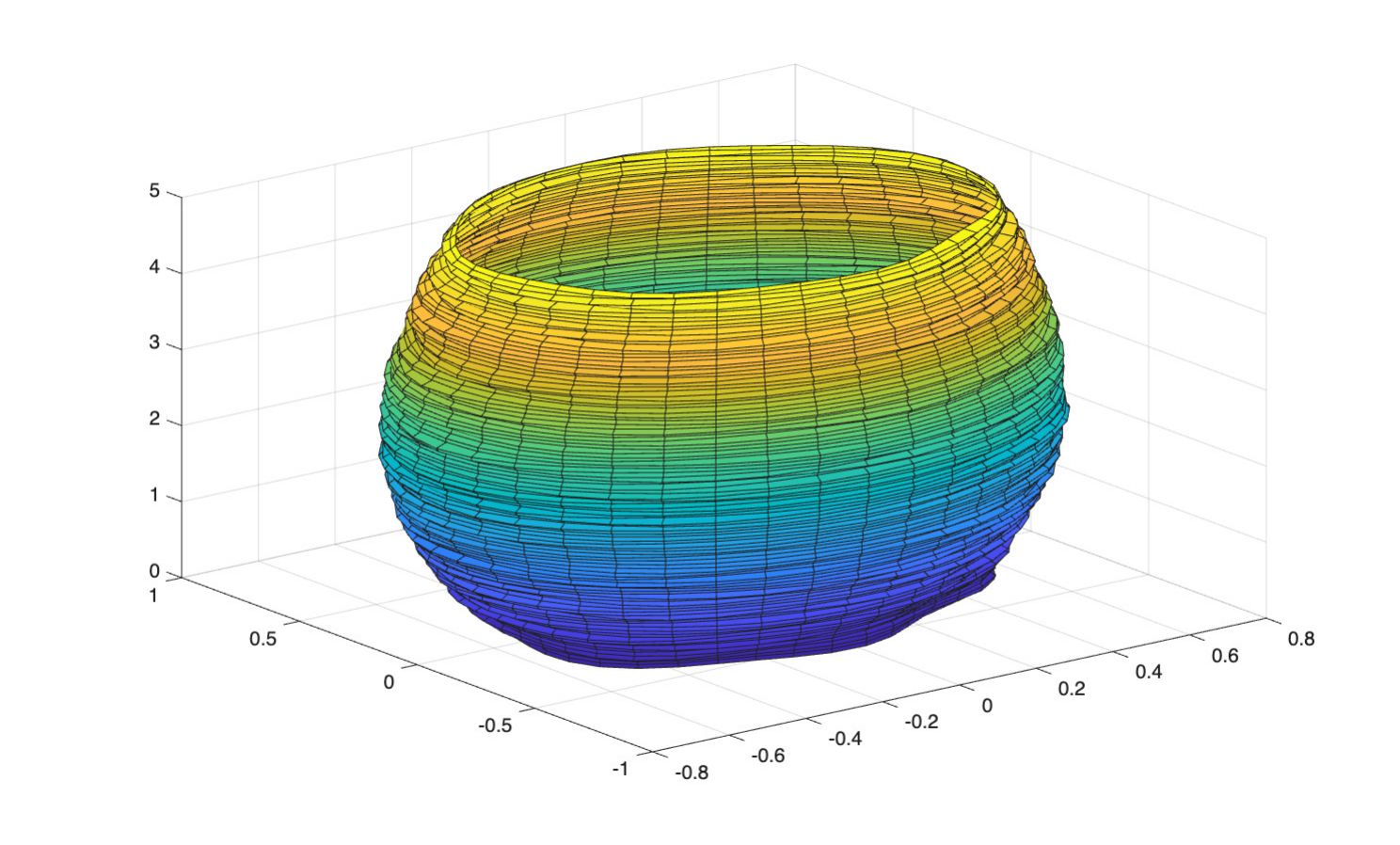}\\
\includegraphics[width=0.48\textwidth,trim={60 30 60 30},clip]{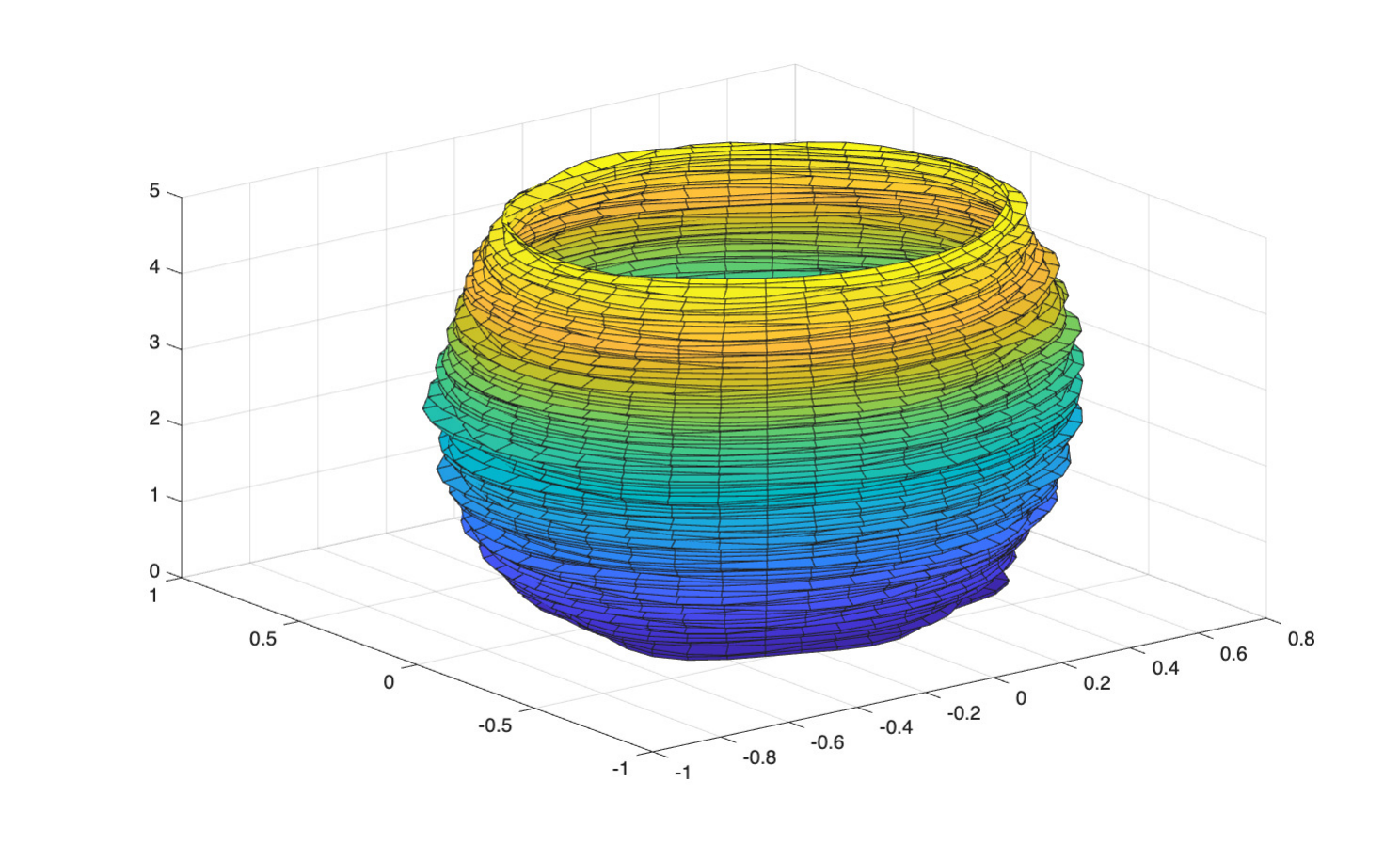}
\includegraphics[width=0.48\textwidth,trim={60 30 60 30},clip]{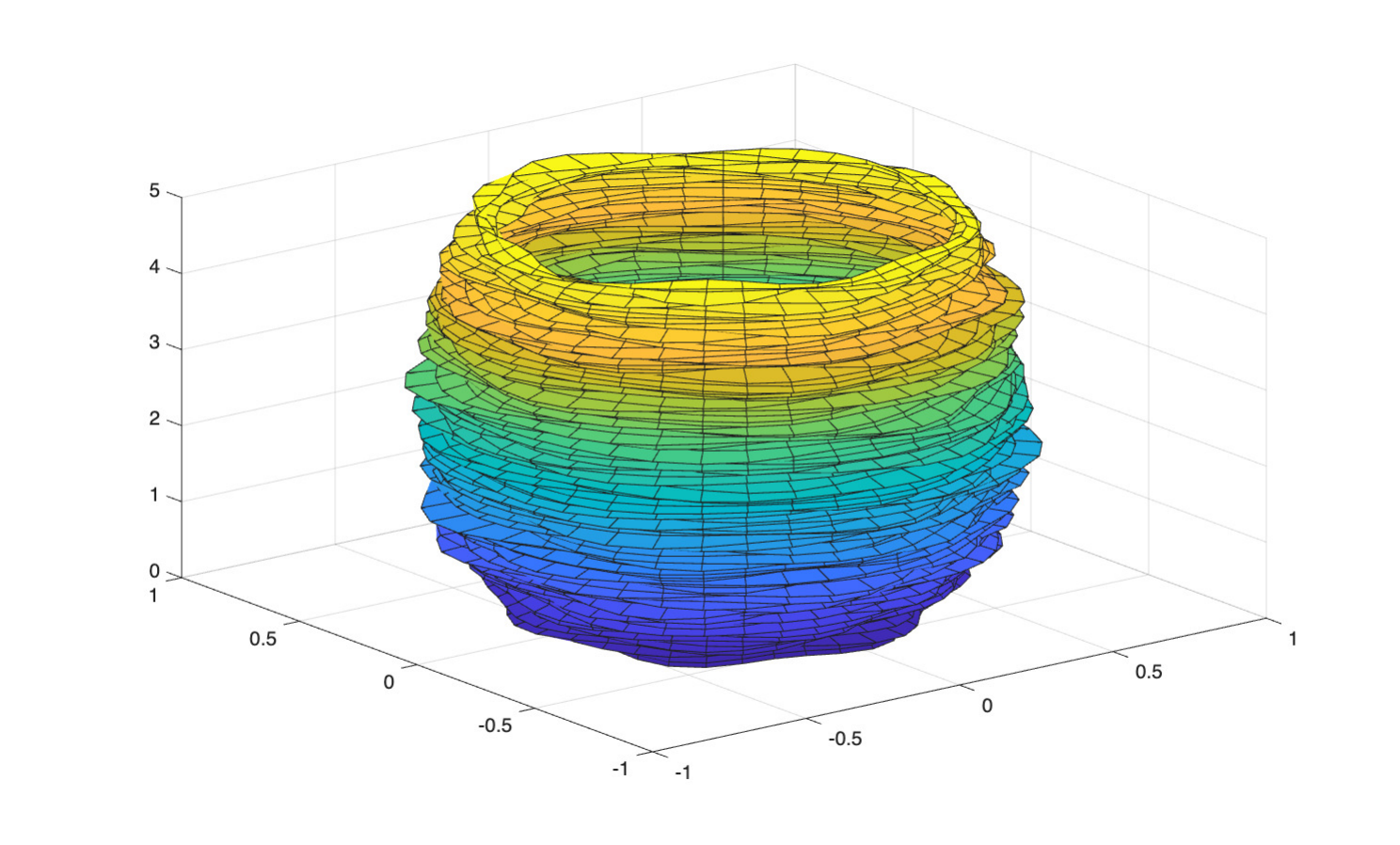}
\caption{\label{fig:inverse1}
The reconstruction of the space-time tube in case of the first example
with $0.25\%$ random noise (top left), $0.5\%$ random noise (top right),
$1.0\%$ random noise (bottom left), and $2.0\%$ random noise (bottom right).}
\end{center}
\end{figure}

\begin{figure}[hbt]
\begin{center}
\includegraphics[width=0.48\textwidth,trim={60 30 60 30},clip]{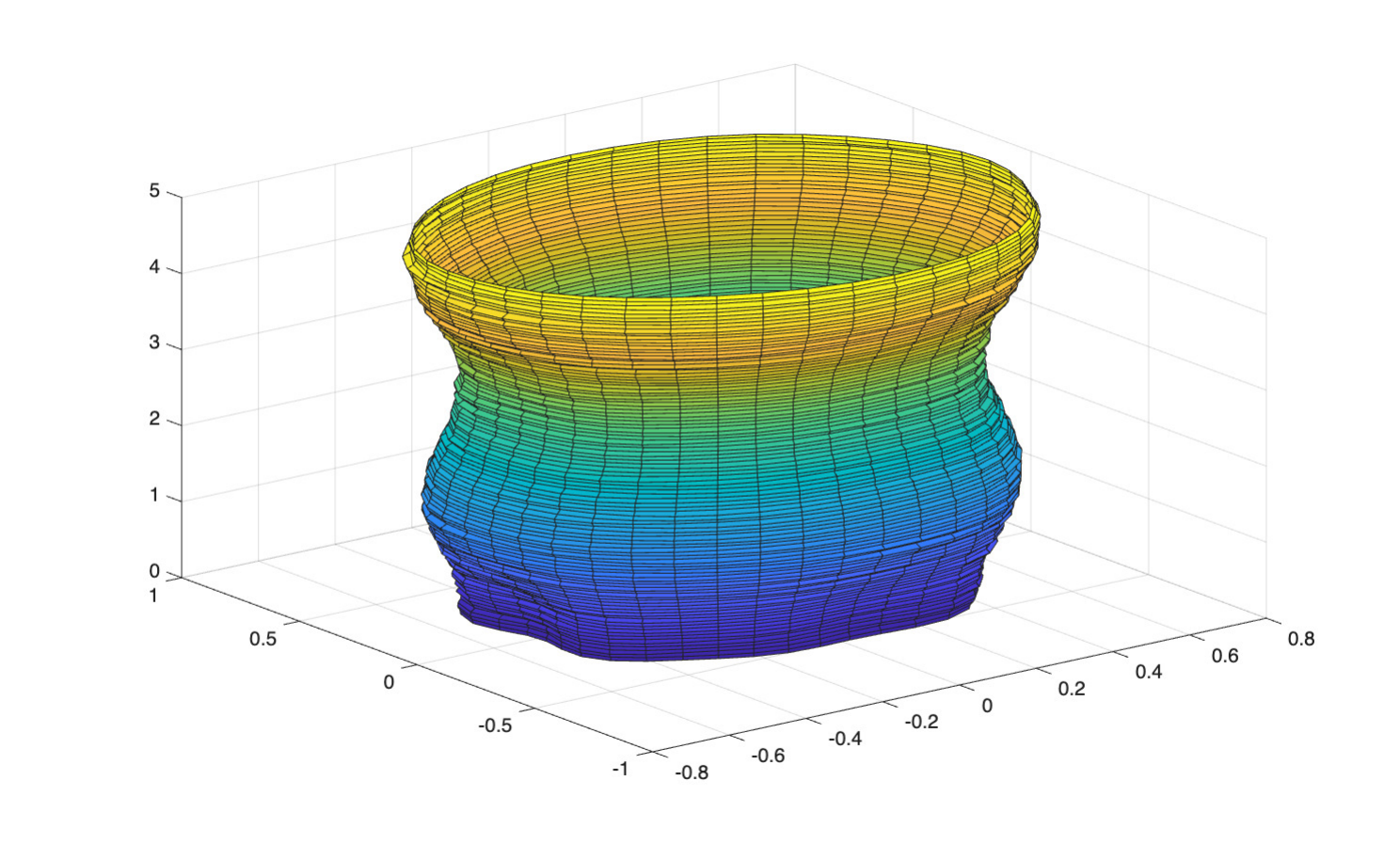}
\includegraphics[width=0.48\textwidth,trim={60 30 60 30},clip]{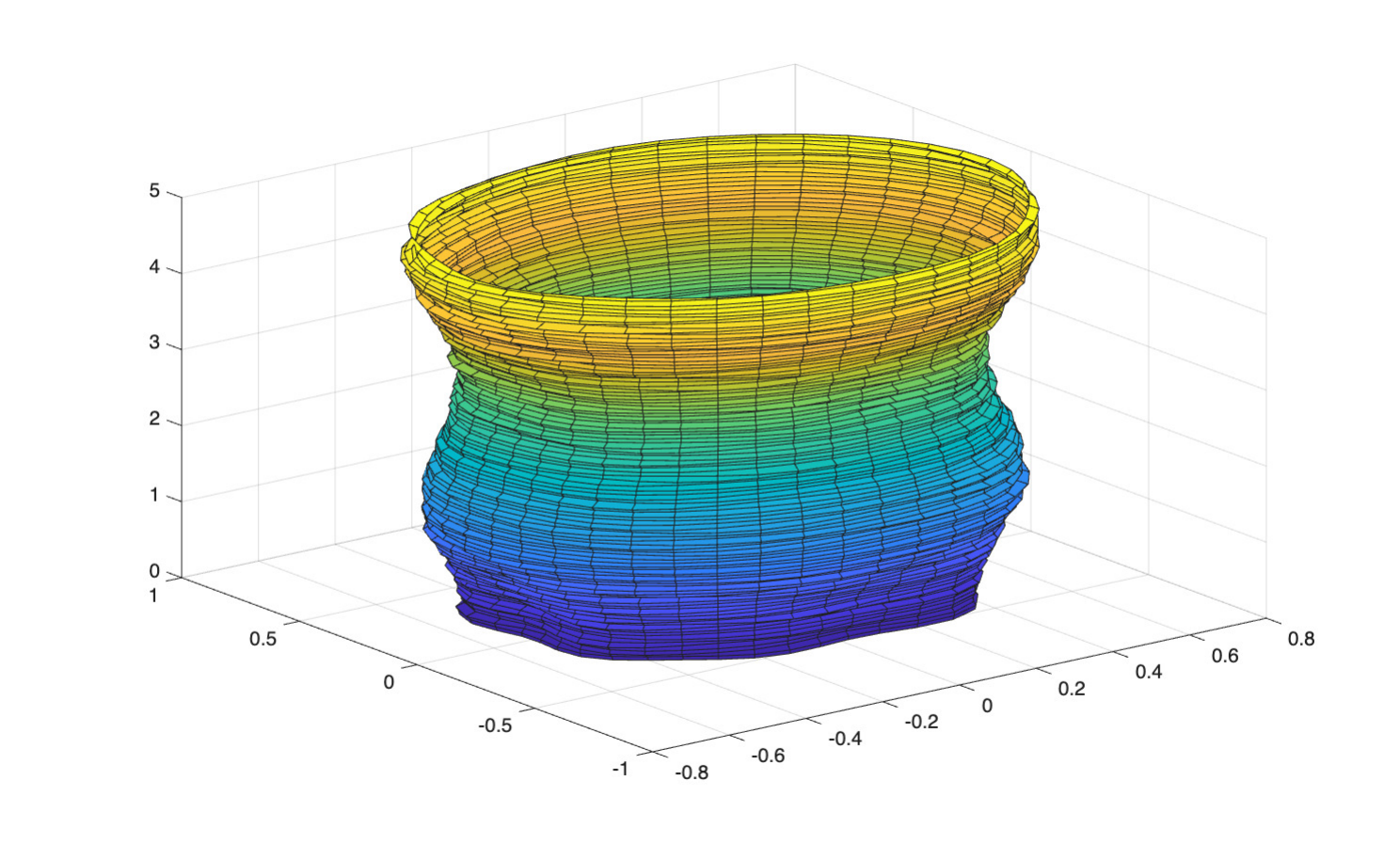}\\
\includegraphics[width=0.48\textwidth,trim={60 30 60 30},clip]{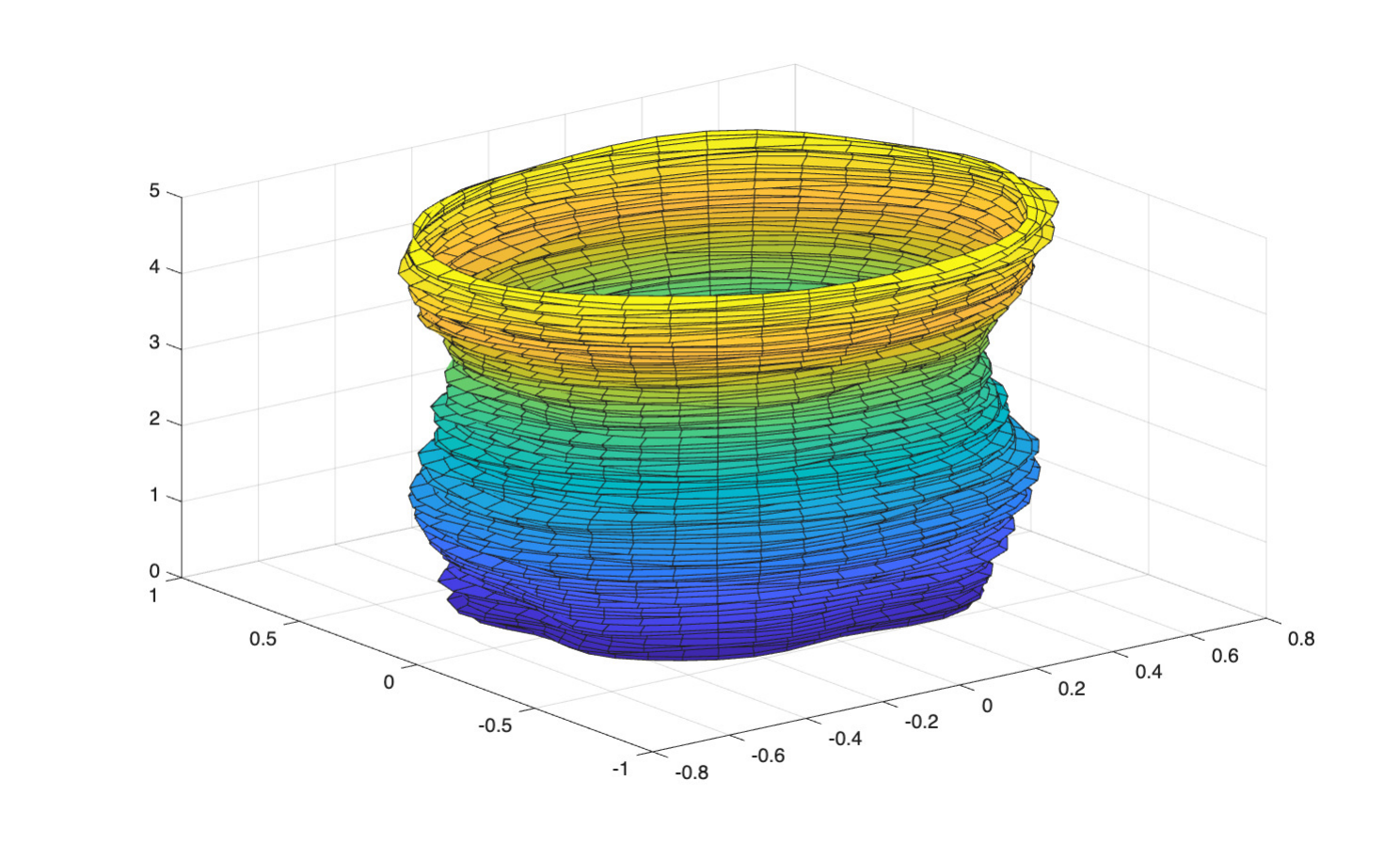}
\includegraphics[width=0.48\textwidth,trim={60 30 60 30},clip]{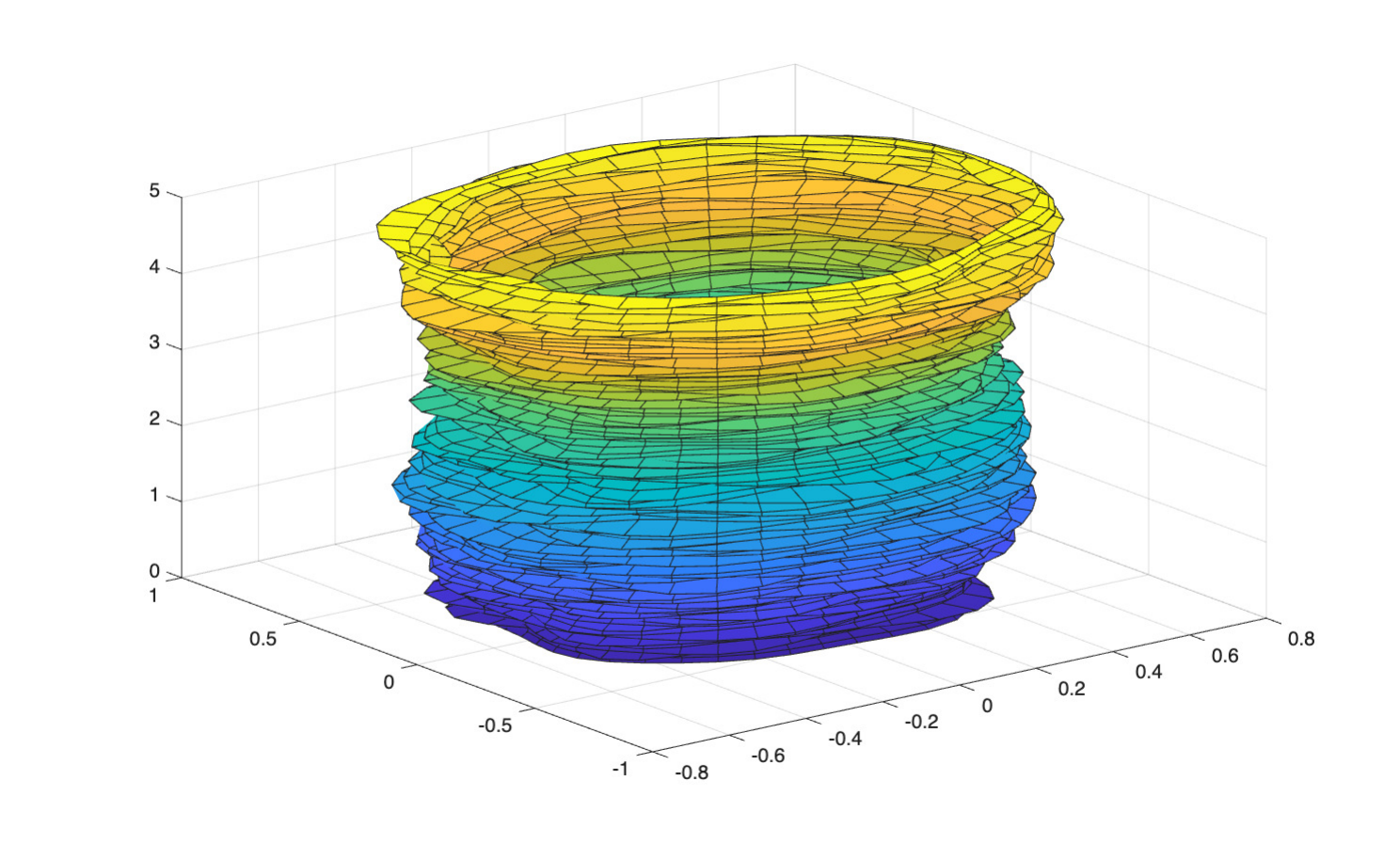}
\caption{\label{fig:inverse2}
The reconstruction of the space-time tube in case of the second example
with $0.25\%$ random noise (top left), $0.5\%$ random noise (top right),
$1.0\%$ random noise (bottom left), and $2.0\%$ random noise (bottom right).}
\end{center}
\end{figure}

The procedure for the reconstruction is as follows.
In the $(k+1)$-st time step, we take the Fourier series 
${\bf a}_k$ for granted to define the domain $\Omega_k$. 
We solve \eqref{eq.Stefan_PDE"}--\eqref{eq.Stefan_initial_cond"}
for the right-hand side $u_m' = 1$ and determine the scaling
factor $\alpha_k$ in such a way that the update \eqref{eq:update}
matches the desired boundary ${\bf a}_{k+1}$:
\begin{equation}\label{eq:reconstruction}
 {\bf x}_{k+1}^{\text{desired}}\overset{!}{=} {\bf x}-\alpha_k\Delta t \frac{\partial v_{k+1}}{\partial{\bf n}}
  	\langle{\bf n},\widehat{\bf x}\rangle\widehat{\bf x}\quad\text{on $\Gamma_k$}.
\end{equation}
Of course, this equation cannot be solved exactly. Hence,
we solve it for all vertices of the finite element mesh that lie
on the boundary in a least-squares sense. 

\subsection{Numerical realization}
In the following, we discuss what the solution of 
\eqref{eq:reconstruction} means in the context of our 
finite element discretization. In view of the linear system 
\eqref{eq:LSE} of equations, we split the the discrete flux 
$\frac{\partial v_{k+1}}{\partial{\bf n}}$ into the contribution 
$\frac{\partial v_{k+1}^{(1)}}{\partial{\bf n}}$ arising from
\[
  \bigg({\bf M}_k+\frac{\Delta t}{2}{\bf A}_k\bigg){\bf v}_{k+1}^{(1)}
  = \bigg({\bf M}_k-\frac{\Delta t}{2}{\bf A}_k\bigg){\bf v}_k,
\]
and the contribution $\frac{\partial v_{k+1}^{(2)}}{\partial{\bf n}}$ 
arising from
\[
  \bigg({\bf M}_k+\frac{\Delta t}{2}{\bf A}_k\bigg){\bf v}_{k+1}^{(2)}
		= -\Delta t {\bf f}_k.
\]
Recall that $\{{\bf y}_{k,i}\}_i$ denote the vertices
of the finite element mesh on $\Omega_k$ and generate
the desired vertices $\{{\bf y}_{k+1,i}^{\text{desired}}\}_i$ 
given by the Fourier coefficients of the measurement 
for $\Omega_{k+1}^{\text{desired}}$. Since they are
equidistant with respect to the polar angle, we can use
the fast Fourier transform for this. Thus, the discrete 
version of \eqref{eq:reconstruction} becomes
\[
 \alpha_k\frac{\partial v_{k+1}^{(2)}}{\partial{\bf n}}({\bf y}_{k,i})
  	\langle{\bf n},{\bf y}_{k,i}\rangle {\bf y}_{k,i}
  = \frac{{\bf y}_{k+1,i}^{\text{desired}}-{\bf y}_{k,i}}{\Delta t}
  -\frac{\partial v_{k+1}^{(1)}}{\partial{\bf n}}({\bf y}_{k,i})
  	\langle{\bf n},{\bf y}_{k,i}\rangle {\bf y}_{k,i}
\]
for $i=1,\ldots,L$ such that the vertex ${\bf y}_{k,i}$ lies
on the boundary $\Gamma_k$. These is an overdetermined
system for the single unknown $\alpha_k$ which we solve
in the least-squares sense.

In Figure~\ref{fig:inverse1} and Figure~\ref{fig:inverse2}, 
we visualized the boundaries fitted in accordance with 
the proposed method for both examples and for all noise 
levels applied. However, we should mention that the
computations in the $(k+1)$-st step is always based 
on the prescribed domain $\Omega_k^{\text{desired}}$.

\subsection{Reconstructions of the melting temperature}
\begin{figure}[hbt]
\includegraphics[width=0.48\textwidth]{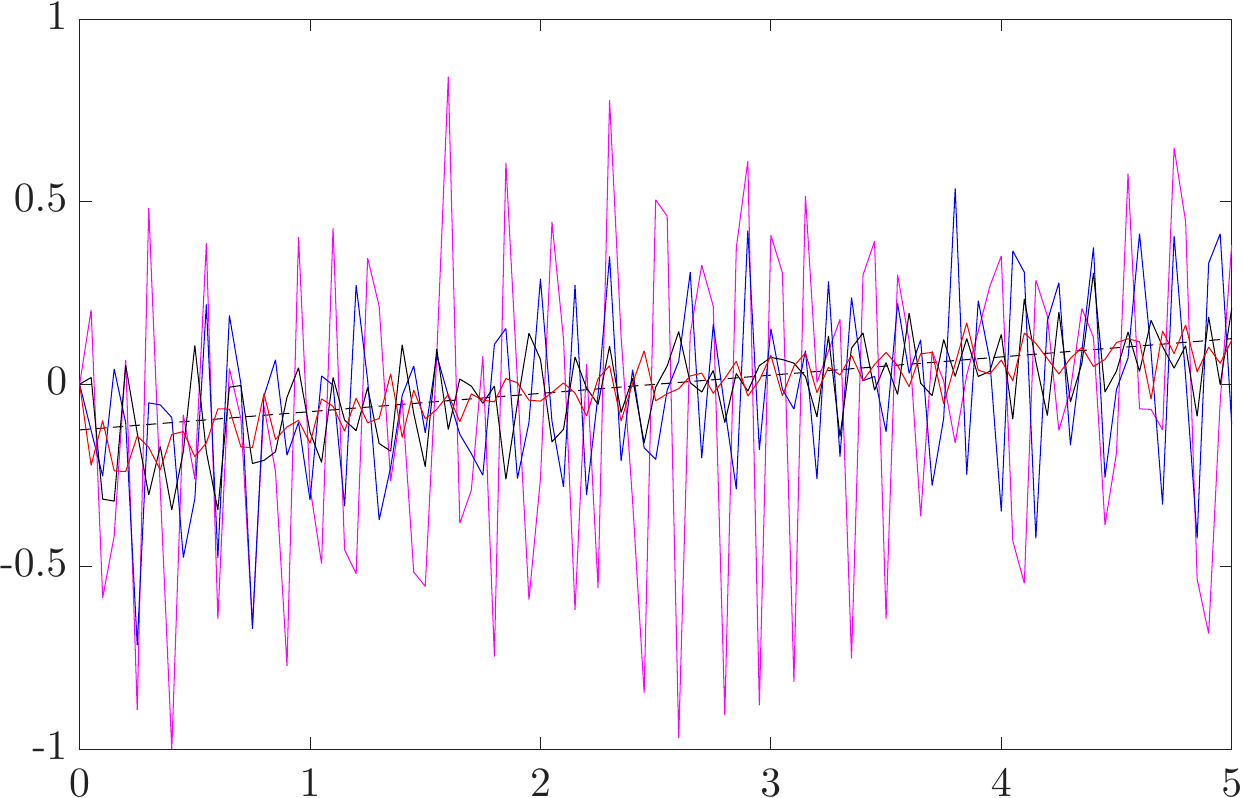}
\includegraphics[width=0.48\textwidth]{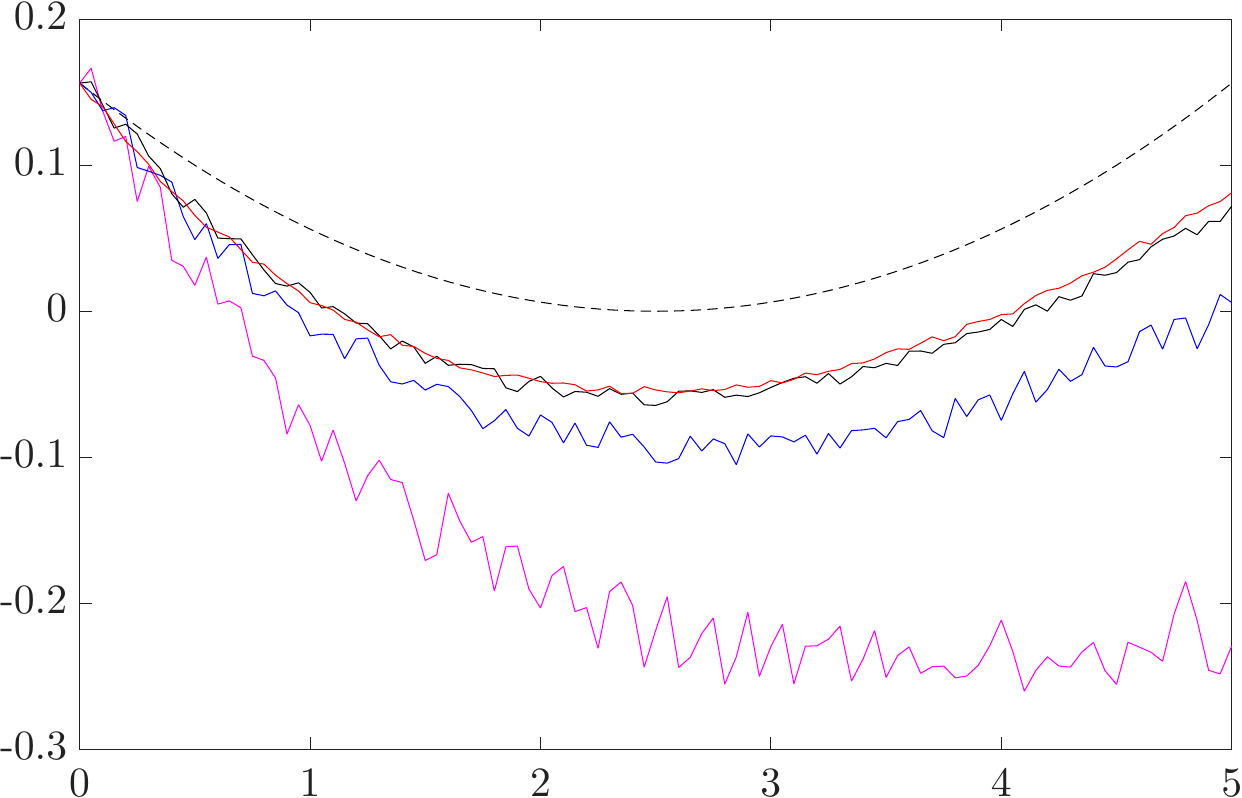}
\caption{\label{fig:melting1}
The reconstruction of the derivative $\dot{u}_m(t)$ of the melting temperature
(left) and the melting temperature $u_m(t)$ (right) in case of the first example
for the different noise levels: $0.25\%$ noise --  red, $0.5\%$ noise --  black, 
$1.0\%$ noise --  blue, and $2.0\%$ noise --  magenta. The true solution is
indicated by the dashed black line.}
\end{figure}

In case of the first example, the reconstructed melting temperature is
found in Figure~\ref{fig:melting1}. The left plot shows the computed
scaling factors which correspond to $\dot{u}_m$. The red curve corresponds
to the noise level $\delta = 0.1\%$, the black one to the noise level $\delta 
= 0.5\%$, the blue one to the noise level $\delta = 1\%$, and the magenta 
one to the noise level $\delta = 2\%$. Especially the latter one oscillates 
quite a lot. However, if we integrate over time to obtain the melting 
temperature $u_m$, the fluctuations are significantly reduced.

We clearly can recognize the applied melting temperature $u_m(t) = 
\nicefrac{1}{20}(t-\nicefrac{5}{2})^2$ in the right plot of Figure~\ref{fig:melting1}.
However, we observe that the discrepancy of the reconstructed melting 
temperature and the prescribed melting temperature grows as the time 
$t$ increases. Moreover, the error of the reconstructions increases as the 
noise level increases as can be expected. Especially, we observe that the 
reconstruction of the melting temperature is not very useful for large time 
$t$ in the case of the highest noise level.

\begin{figure}[hbt]
\includegraphics[width=0.48\textwidth]{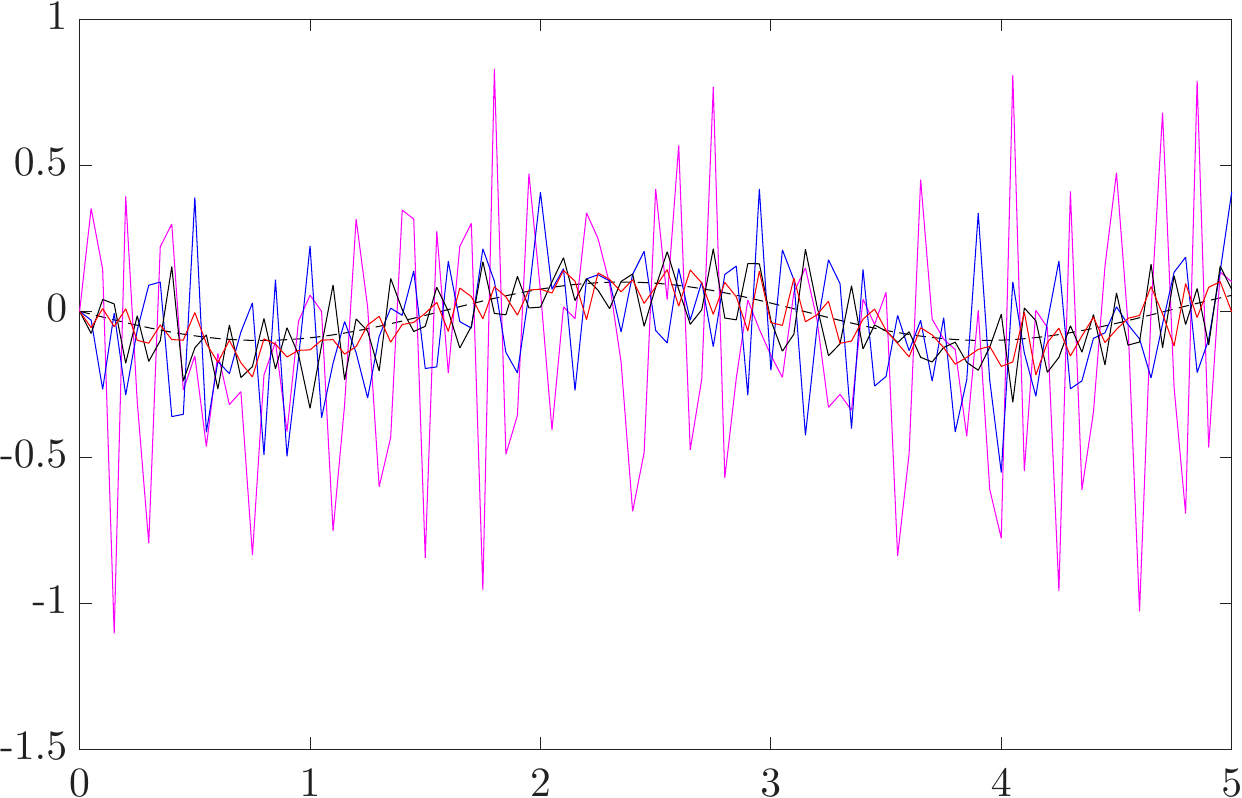}
\includegraphics[width=0.48\textwidth]{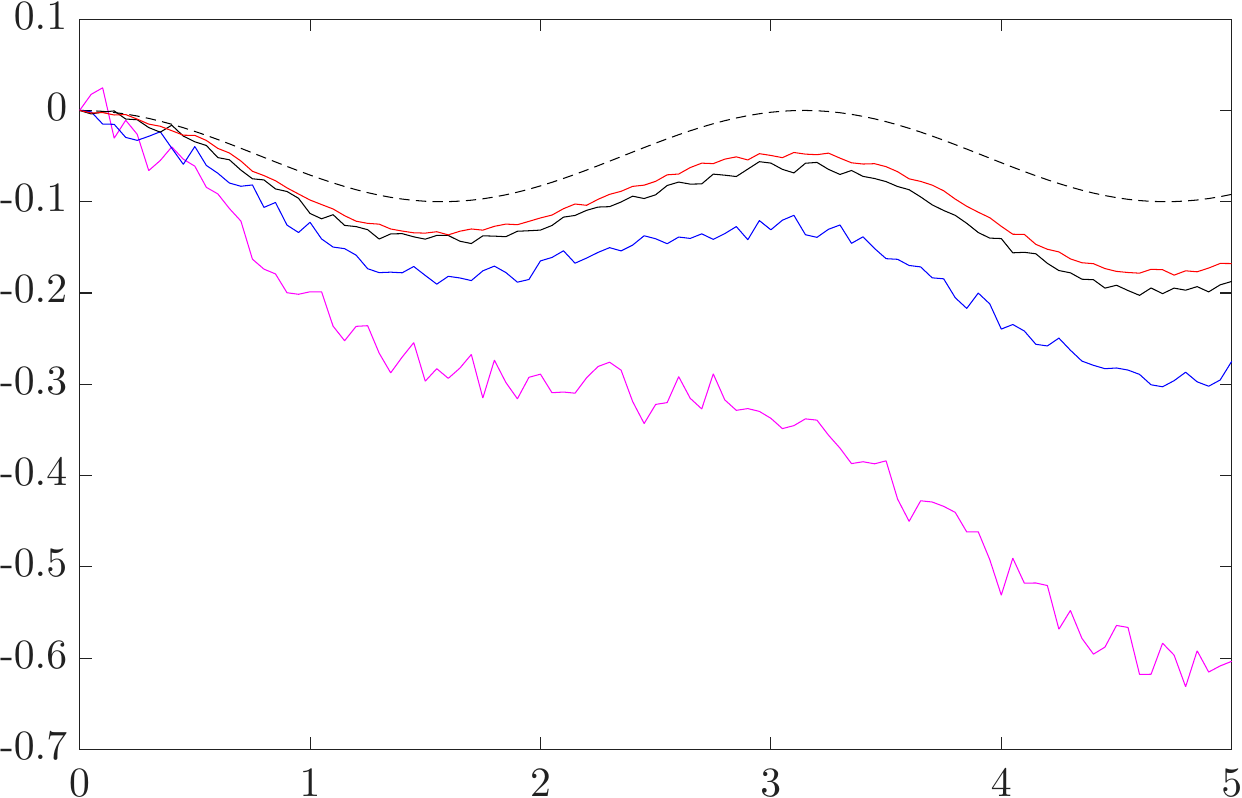}
\caption{\label{fig:melting2}
The reconstruction of the derivative $\dot{u}_m(t)$ of the melting temperature
(left) and the melting temperature $u_m(t)$ (right) in case of the second example
for the different noise levels: $0.25\%$ noise --  red, $0.5\%$ noise --  black, 
$1.0\%$ noise --  blue, and $2.0\%$ noise --  magenta. The true solution is
indicated by the dashed black line.}
\end{figure}

Similar observations can be made in case of the reconstructions 
of the melting temperature in the second example, which are found 
in Figure~\ref{fig:melting2}. We likewise observe in the left plot of
Figure~\ref{fig:melting2} an increase of the oscillations in $\dot{u}_m$
when the noise level increases, where the colour coding of the 
graphs is the same as for the first example. Nonetheless, the 
melting temperature $u_m(t)$ itself is reasonably well reconstructed 
as can be seen in the right plot of Figure~\ref{fig:melting2}. But the 
reconstruction of the melting temperature for the highest noise 
level is again very inaccurate for larger time $t$.

\section{Conclusion}
\label{sct:conclusio} 
In the present article, we considered the numerical solution 
of a Stefan problem, where the melting temperature varies
over time. This happens in practical applications when the 
pressure in the external space changes during time for 
example. We provided a numerical algorithm to simulate 
such a transient Stefan problem that is based on finite 
elements defined on a moving mesh. We investigated 
in addition also the respective inverse problem -- the 
reconstruction of the melting temperature from measurements 
of the evolving boundary. We demonstrated that this inverse 
problem can indeed be solved numerically, leading to 
reasonable reconstructions of the melting temperature. 
We presented different computations to validate the 
feasibility of the proposed approach.

\bibliographystyle{plain}
\bibliography{literature}
\end{document}